\documentclass[12pt,reqno]{amsart}

\usepackage{mathabx}
\usepackage{amssymb}
\usepackage{bbm}
\usepackage{mathrsfs}
\usepackage[all]{xy}
\usepackage{float}
\usepackage{graphicx}

\newtheorem{theo}{Theorem}[section]
\newtheorem{lem}[theo]{Lemma}

\newtheorem{coro}[theo]{Corollary}
\newtheorem{conj}[theo]{Conjecture}

\theoremstyle{definition}
\newtheorem{defi}[theo]{Definition}

\newtheorem{ttt}[theo]{}

\theoremstyle{remark}
\newtheorem{rem}[theo]{Remark}
\newtheorem{rems}[theo]{Remarks}
\newtheorem{nota}[theo]{Notation}
\newtheorem{ex}[theo]{Example}
\newtheorem{exs}[theo]{Examples}

\newcommand{\br}{ }
\newcommand{\brr}{, }

\newcommand{\C}{\mathop{\text{\rm C}}\nolimits}
\newcommand{\E}{\mathop{\text{\rm E}}\nolimits}
\newcommand{\K}{\mathop{\text{\rm K}}\nolimits}
\newcommand{\Ima}{\mathop{\text{\rm Im}}\nolimits}
\newcommand{\Gal}{\mathop{\text{\rm Gal}}\nolimits}
\newcommand{\Cl}{\mathop{\text{\rm Cl}}\nolimits}
\newcommand{\Spec}{\mathop{\text{\rm Spec}}\nolimits}
\newcommand{\Frob}{\mathop{\text{\rm Frob}}\nolimits}
\newcommand{\Pic}{\mathop{\text{\rm Pic}}\nolimits}
\newcommand{\GL}{\mathop{\text{\rm GL}}\nolimits}

\renewcommand{\O}{\mathop{\text{\rm O}}\nolimits}
\newcommand{\G}{\mathop{\text{\rm G}}\nolimits}
\newcommand{\GO}{\mathop{\text{\rm GO}}\nolimits}
\newcommand{\SO}{\mathop{\text{\rm SO}}\nolimits}
\newcommand{\U}{\mathop{\text{\rm U}}\nolimits}

\newcommand{\ST}{\mathop{\text{\rm ST}}\nolimits}
\newcommand{\Ad}{\mathop{\text{\rm Ad}}\nolimits}
\newcommand{\Hg}{\mathop{\text{\rm Hg}}\nolimits}
\newcommand{\Lie}{\mathop{\text{\rm Lie}}\nolimits}
\newcommand{\im}{\mathop{\text{\rm im}}\nolimits}
\newcommand{\rk}{\mathop{\text{\rm rk}}\nolimits}
\renewcommand{\Re}{\mathop{\text{\rm Re}}\nolimits}
\newcommand{\den}{\mathop{\text{\rm den}}\nolimits}

\newcommand{\disc}{\mathop{\text{\rm disc}}\nolimits}

\newcommand{\Zar}{\text{\rm Zar}}
\newcommand{\Haar}{\text{\rm Haar}}
\newcommand{\un}{\text{\rm un}}
\newcommand{\sss}{\text{\rm ss}}
\newcommand{\alg}{\text{\rm alg}}
\newcommand{\tr}{\text{\rm tr}}
\newcommand{\Tr}{\text{\rm Tr}}
\newcommand{\et}{\text{\rm \'et}}
\newcommand{\id}{\text{\rm id}}
\newcommand{\pr}{\text{\rm pr}}
\newcommand{\bl}{\text{\rm bl}}
\newcommand{\diag}{\text{\rm diag}}

\newcommand{\triv}{\text{\rm triv}}

\newcommand{\Pb}{{\text{\bf P}}}
\newcommand{\Vb}{{\text{\bf V}}}

\newcommand{\bbC}{{\mathbbm C}}
\newcommand{\bbF}{{\mathbbm F}}
\newcommand{\bbN}{{\mathbbm N}}
\newcommand{\bbP}{{\mathbbm P}}
\newcommand{\bbQ}{{\mathbbm Q}}
\newcommand{\bbR}{{\mathbbm R}}
\newcommand{\bbZ}{{\mathbbm Z}}

\newcommand{\calL}{{\mathscr{L}}}
\newcommand{\calO}{{\mathscr{O}}}
\newcommand{\calX}{{\mathscr{X}}}

\newcommand{\frakn}{\mathfrak{n}}

\newcommand{\Exterior}{\mathchoice{{\textstyle\bigwedge}}%
 {{\bigwedge}}%
 {{\textstyle\wedge}}%
 {{\scriptstyle\wedge}}}

\newcounter{abc}
\newenvironment{abc}{\begin{list}{\rm \alph{abc}) }%
{\usecounter{abc} \leftmargin=0.0pt \labelsep=0.0pt %
\listparindent=0.0pt \labelwidth=0.0pt \parsep=\smallskipamount%
 \itemsep=0.0pt \topsep=0.0pt \partopsep=\smallskipamount}}{\end{list}}

\newcounter{iii}
\newenvironment{iii}{\begin{list}{\rm \roman{iii}) }%
{\usecounter{iii} \leftmargin=0.0pt \labelsep=0.0pt %
\listparindent=0.0pt \labelwidth=0.0pt \parsep=\smallskipamount%
 \itemsep=0.0pt \topsep=0.0pt \partopsep=\smallskipamount}}{\end{list}}

\def\rightend#1#2{{%
 \leavevmode\nobreak\hskip .5em plus 1fil
 \penalty600 \hskip 0pt plus -1filll
 \vadjust{}\nobreak\hskip 0pt plus 1filll%
 #1\parfillskip=#2\relax \par}}

\def\eop{\ifmmode\rule[-22pt]{0pt}{1pt}\ifinner\tag*{$\square$}\else\eqno{\square}\fi\else\rightend{$\square$}{0pt}\fi
}

\thanks{}

\title[Frobenius trace distributions for $K3$~surfaces]{Frobenius trace distributions for {\boldmath$K3$}~surfaces}

\begin{document}

\author{Andreas-Stephan Elsenhans}

\address{Institut f\"ur Mathematik\\ Universit\"at W\"urzburg\\ Emil-Fischer-Stra\ss e 30\\ D-97074 W\"urzburg\\ Germany}
\email{stephan.elsenhans@mathematik.uni-wuerzburg.de}
\urladdr{https://www.mathematik.uni-wuerzburg.de/institut/personal/elsenhans.html}

\author[J\"org Jahnel]{J\"org Jahnel}

\address{\mbox{Department Mathematik\\ \!Univ.\ \!Siegen\\ \!Walter-Flex-Str.\ \!3\\ \!D-57068 \!Siegen\\ \!Germany}}
\email{jahnel@mathematik.uni-siegen.de}
\urladdr{https://www.uni-math.gwdg.de/jahnel}

\thanks{}

\date{September~21, 2022}



\begin{abstract}
We study the distribution of the Frobenius traces on
$K3$~surfaces.
We~compare experimental data with the predictions made by the Sato--Tate conjecture, i.e.\ with the theoretical distributions derived from the theory of Lie groups assuming equidistribution. Our sample consists of generic
$K3$~surfaces,
as well as of such having real and complex multiplication. Each time, the theoretical density and the histogram obtained by counting points match in the range of visible accuracy. Thus, we report evidence for the Sato--Tate conjecture for the surfaces considered.
\end{abstract}

\maketitle
\thispagestyle{empty}

\section{Introduction}

Given a smooth, projective variety
$X$
over~$\bbQ$,
one may choose a model
$\calX$
of~$X$
that is projective
over~$\Spec\bbZ$.
The~point counts
$\#\calX_p(\bbF_{\!p})$,
at least for the
primes~$p$
of good reduction, then form a highly interesting set of quantities related to the
variety~$X$.
For~example, for
$X$
an elliptic curve, Hasse's bound states that
$a_p \in [-2,2]$,
for
$\smash{a_p := (\#\calX_p(\bbF_{\!p}) - p - 1)/\sqrt{p\mathstrut}}$,
and it seems natural to ask for the distribution of the sequence
$(a_p)_{p\in\bbP}$
in that~interval.\vspace{-2mm}

\begin{figure}[H]
\label{ell_curve_plots}
\begin{center}
\includegraphics[scale=0.4]{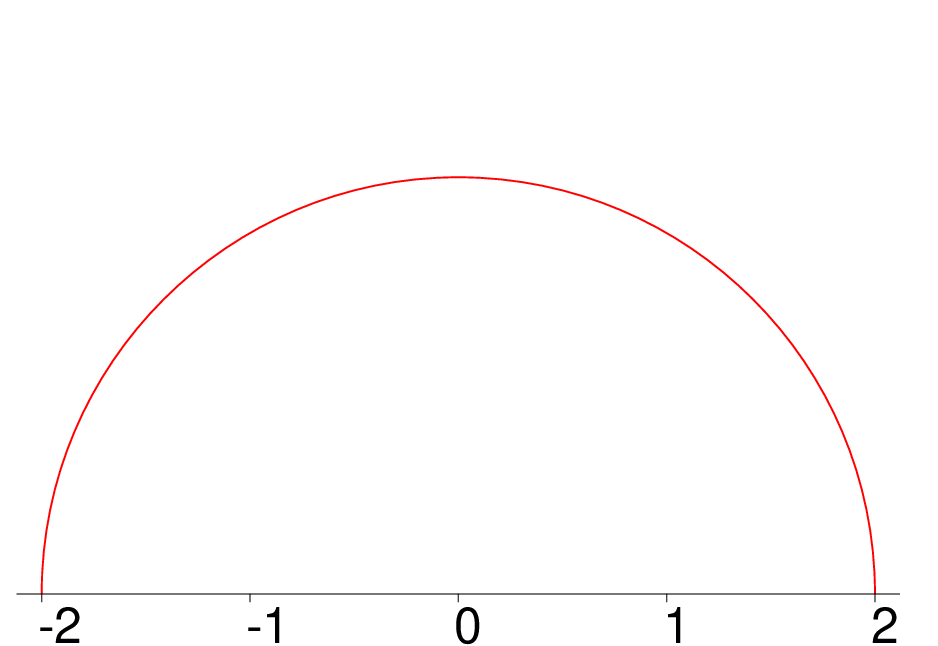}
\hspace{2cm}
\includegraphics[scale=0.4]{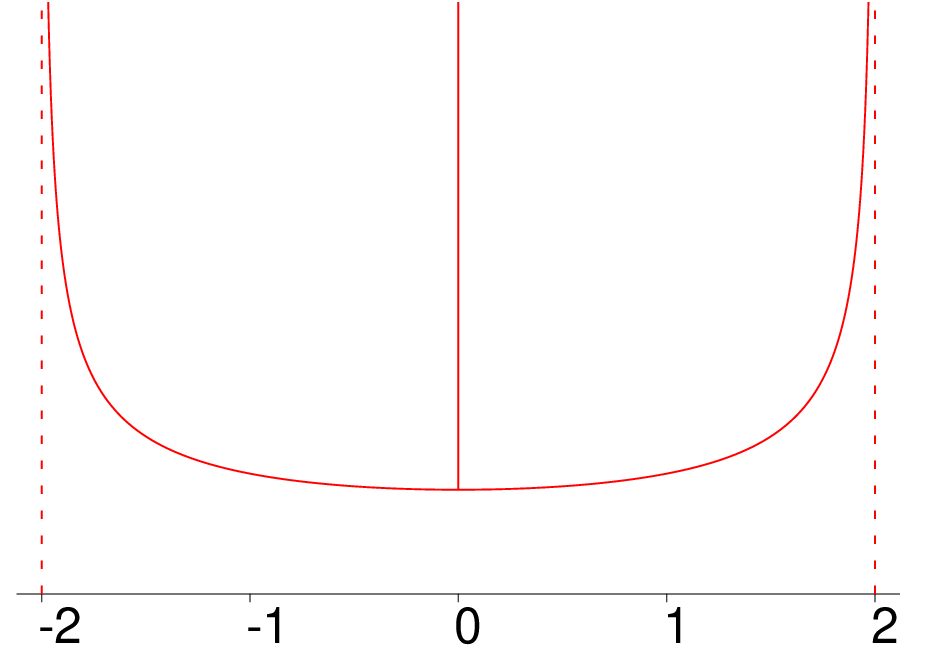}
\end{center}
\caption{Distribution for elliptic curves, general (left) and CM (right)}
\end{figure}\vspace{-5mm}

When~$X$
does not have CM,
$(a_p)_{p\in\bbP}$
is equidistributed with respect to the measure with density
$\smash{\frac1{2\pi}\sqrt{4-t^2}}$.
This~was first observed experimentally by M.\ Sato. J.~Tate gave a partial explanation based on what is now called the Tate conjecture, applied to the direct powers
$X^m$
\mbox{\cite[\S4]{Tat}}. Making~strong assumptions, related to modularity, on the elliptic curve, J.-P.\ Serre \cite[\S I.A.2, Example~3]{Se68} provided a proof shortly afterwards. Based on new developments originating with A.\ Wiles~\cite{Wi}, the result was finally established unconditionally by R.\ Taylor et al.\ \cite{CHT,HST,Tay}.
The CM case is substantially easier and was well understood already in the sixties. Here,~the density function
is~$\smash{\frac1{\pi\sqrt{4-t^2}}}$.
A~proof may be found in~\cite[Proposition~2.16]{Su19}.

For curves of higher genus, extensive experiments have been carried out by A.\ Sutherland~\cite{Su19,KS}. Numerical data are available from his website~\cite{Su20}. Theoretical investigations concerning the
genus-$2$
case are made in~\cite{FKRS}. An~equidistribution statement is proven in certain cases, in which the Jacobian is geometrically isogenous to a direct product of elliptic curves, cf.\ \cite{FS}. For an arbitrary smooth, projective
variety~$X$,
there is the Sato--Tate conjecture, which we describe in Section~\ref{sec_ST}. It~is not just concerned with the traces, but predicts equidistribution of certain
elements~$x_p$
derived from the Frobenii within a compact Lie group, the {\em Sato--Tate group\/}
$\ST^i(X)$.

In~the situation of a
$K3$~surface,
there is an explicit description of the neutral component
$(\ST^2(X))^0$
of the Sato--Tate group, due to the work of Yu.\,G.\ Zarhin~\cite{Za} and S.\,G.\ Tankeev \cite{Tan90,Tan95}. In~fact,
$(\ST^2(X))^0$
depends only on the geometric Picard rank, the degree of the endomorphism
field~$E$,
and the bit of information whether
$E$
is totally real or a CM~field. We~recall the description of
$(\ST^2(X))^0$
in Corollary~\ref{ST0}.

Note~that the concept of the endomorphism field is more subtle here than for abelian varieties, since one considers the endomorphism field of a Hodge structure associated
with~$X$.
Cf.~Paragraph~\ref{neut_comp} for details.

\subsubsection*{A theoretical result}
For~arbitrary
$K3$~surfaces,
we give an upper bound for the possible component groups
of~$\ST^2(X)$
in Theorem~\ref{upper_bound}. For example, in the case of real multiplication by a quadratic number field,
$\smash{\ST^2_\tr(X)/(\ST^2(X))^0}$
is naturally contained in the dihedral group of order eight. We show that the order is, in fact, at most four and indicate that two {\em non-isomorphic\/} subgroups of order four are~possible.\medskip

But~this is not the main goal of this~article. The~main goal  is to report on our experiments concerning the Sato--Tate conjecture for certain
$K3$~surfaces
of geometric Picard rank
$16$
(and~$17$).
The~surfaces in our sample have singular models of degree two and vary in endomorphism field and jump character~\cite{CEJ}. For every surface, we determine the Sato--Tate measure, i.e.\ the theoretical distribution of the Frobenius traces, according to the Sato--Tate conjecture, and compare it with a histogram obtained by explicitly counting points, for all
primes~$p$
up
to~$10^8$.

\subsubsection*{The experimental results}
For any of the seven surfaces in our sample, the theoretical distribution and the histogram match in the range of visible accuracy. We~present the seven histograms, each one in juxtaposition with the graph of the corresponding theoretical density function, in Section~\ref{sec_exp}. Cf.\ Figures~2 to~5. The~same information in higher resolution is available from the second author's web page at {\tt https://www.uni-math.gwdg.de/jahnel/Arbeiten/histograms.tar.gz}.

Concerning the rate of convergence, our data suggest that the order
is~$\smash{\frac12}$.
I.e., that the distribution of the Frobenius traces converges towards the Sato--Tate measure of
order~$\smash{\frac12}$,
in terms of the number of primes~used. We~present data supporting such a conjecture for one of our examples. The~other surfaces show qualitatively the same~behaviour.

\subsubsection*{The selection of our sample}
Any~selection of examples is, of course, somewhat~arbitrary. However,~we strongly feel that
$K3$~surfaces
of geometric Picard
rank~$16$
are a very reasonable compromise between surfaces of high rank, which are very special, and surfaces of low rank, which are certainly general, but hard to~treat. Most~notably,
rank~$16$
is the largest one that allows real multiplication~\cite[Lemma~3.2]{vG}.

We~include an example of geometric Picard
rank~$17$
and trivial jump character, as, in this situation, the theoretical density function is not symmetric. Cf.\ the third histogram in Figure~2.

\subsubsection*{Computations}

The computations related to this project were done with {\tt magma} \cite{BCP}. For~symbolic integration, we used {\tt maple}~\cite{Ma}.

\section{The Sato--Tate conjecture}
\label{sec_ST}

\subsubsection*{The algebraic monodromy group}

Let
$X$
be a smooth, projective variety
over~$\bbQ$,
and
$S$
a finite set of primes, outside of which
$X$
has good~reduction. Then~the Lefschetz trace formula in \'etale cohomology~\cite[Expos\'e~III, Th\'eor\`eme 6.13.3]{SGA5}, together with the smooth specialisation theorem~\cite[Expos\'e~XVI, Corollaire~2.2]{SGA4}, show
\begin{equation}
\label{lefschetz}
\#\calX_p(\bbF_{\!p}) = \!\sum_{i=0}^{2\dim X}\!\! (-1)^i \,\Tr \big(\! \Frob_p\colon H_\et^i(X_{\overline\bbQ}, \bbQ_l) \righttoleftarrow\!\big) \,,
\end{equation}
for
$p \in \bbP \setminus S$
and any
prime~$l \neq p$.
Here,
$\Frob_p \in \Gal(\overline\bbQ/\bbQ)$
denotes a Frobenius~lift.
Since~$\Frob_p$
is unique up to conjugation, the trace is independent of that~choice.

Suppose now that
$i$
is even, which is the slightly easier case and the one we study in this~article. Then
\begin{equation}
\label{tate}
\Tr \big(\! \Frob_p\colon H_\et^i(X_{\overline\bbQ}, \bbQ_l) \righttoleftarrow\!\big) = p^{\frac{i}2} \!\cdot \Tr \big(\! \Frob_p\colon H_\et^i(X_{\overline\bbQ}, \bbQ_l(i/2)) \righttoleftarrow\!\big) \,.
\end{equation}
The trace on the right hand side is known to be a rational number that is independent
of~$l$.
According to the Weil conjectures, proven by P.\ Deligne~\cite[Th\'eor\`eme~1.6]{De74}, every eigenvalue
of~$\Frob_p$
is an algebraic number, all complex (and real) embeddings of which are of absolute
value~$1$.

Moreover,~the operation of the Galois~group,
\begin{equation}
\label{gal_rep}
\varrho^i_{X,l}\colon \Gal(\overline\bbQ/\bbQ) \longrightarrow \GL(H_\et^i(X_{\overline\bbQ}, \bbQ_l(i/2))) \,,
\end{equation}
is continuous, cf.~\cite[Expos\'e~VIII, Th\'eor\`eme~5.2]{SGA4}. Its image is hence an
\mbox{$l$-adic}
Lie~group. The Zariski closure
$\smash{G_{X,l}^{i,\Zar} := \overline{\im(\varrho^i_{X,l})}}$
is called the {\em algebraic monodromy group\/}
of~$X$
(in
degree~$i$).
It~is a linear algebraic group
over~$\bbQ_l$.
%
%The component group
%$\smash{G_{X,l}^{i,\Zar}/(G_{X,l}^{i,\Zar})^0}$,
%which is always finite, is known to be independent
%of~$l$~\mbox{\cite[p.~15ff]{Se81}}, cf.\ \cite[\S8.3.4]{Se12}.

\subsubsection*{Inclusion in the orthogonal group}

Fix a hyperplane section
$H \subset X$.
Then, by Poincar\'e duality and the hard Lefschetz theorem \cite[Th\'eor\`eme~4.1.1]{De80}, the cohomology vector space
$H_\et^i(X_{\overline\bbQ}, \bbQ_l(i/2))$
is equipped with a non-degenerate, symmetric, bilinear pairing. For
$i \leq \dim X$,
this is given as~follows,
\begin{eqnarray*}
H_\et^i(X_{\overline\bbQ}, \bbQ_l(i/2)) \times H_\et^i(X_{\overline\bbQ}, \bbQ_l(i/2)) & \longrightarrow & H_\et^{2\dim X}(X_{\overline\bbQ}, \bbQ_l(\dim X)) \cong \bbQ_l \,, \\
(\alpha, \beta) & \mapsto & \langle\alpha, \beta\rangle := \alpha \cup \beta \cup [H]^{\dim X - i} \,.
\end{eqnarray*}
The operation of
$\Gal(\overline\bbQ/\bbQ)$
respects this pairing, so one actually has an inclusion
$$G_{X,l}^{i,\Zar} \subseteq \O(H_\et^i(X_{\overline\bbQ}, \bbQ_l(i/2))) \,.$$

\subsubsection*{The Sato--Tate group}

Let us fix an embedding
$\bbQ_l \hookrightarrow \bbC$.
Then
$\smash{G_{X,l}^{i,\Zar}(\bbC)}$
is a complex Lie group, equipped with an inclusion
$\smash{\iota\colon G_{X,l}^{i,\Zar} \hookrightarrow G_{X,l}^{i,\Zar}(\bbC)}$.
Moreover,
$\smash{G_{X,l}^{i,\Zar}(\bbC)}$
is contained in the matrix group
$\O(H_\et^i(X_{\overline\bbQ}, \bbQ_l(i/2)) \otimes_{\bbQ_l}\! \bbC)$.
In particular, the
elements of~$\smash{G_{X,l}^{i,\Zar}(\bbC)}$
have eigenvalues being complex numbers and there is the trace~map
$$\tr\colon G_{X,l}^{i,\Zar}(\bbC) \to \bbC \,.$$

The~maximal compact subgroup
$\ST^i(X)$
of~$\smash{G_{X,l}^{i,\Zar}(\bbC)}$
is called the {\em Sato--Tate group\/}
of~$X$
in
degree~$i$.
The Sato--Tate group is a compact Lie group, in general disconnected. For~the component group, one clearly~has
$$\ST^i(X)/(\ST^i(X))^0 \cong G_{X,l}^{i,\Zar}(\bbC)/(G_{X,l}^{i,\Zar}(\bbC))^0 \cong G_{X,l}^{i,\Zar}/(G_{X,l}^{i,\Zar})^0 \,.$$

\begin{rems}
\begin{iii}
\item
The maximal compact subgroups of a Lie group with finitely many connected components are mutually conjugate~\cite[Theorem IV.3.5]{OV}. Thus, the Sato--Tate group is well-defined, up to~conjugation.
\item
According to the Mumford--Tate conjecture, the neutral component
$\smash{(G_{X,l}^{i,\Zar})^0}$
of the algebraic monodromy group
$\smash{G_{X,l}^{i,\Zar}}$
coincides with
$\Hg^i(X) \!\times_{\Spec\bbQ}\! \Spec\bbQ_l$,
for~$\Hg^i(X)$
the
$i$-the
Hodge group
of~$X$
\cite[Definition\ 3.8 and Conjecture~3.10]{Su19}.
\end{iii}
\end{rems}

\begin{rem}
One might want to work without Tate twist, as one is forced to do in the case when
$i$
is~odd. The algebraic monodromy group is then only contained in
$\GO(H_\et^i(X_{\overline\bbQ}, \bbQ_l))$
and one would impose an orthogonal constraint, i.e.\ intersect with the orthogonal group, afterwards. Cf.~\cite{Su19}, in particular \cite[Remark~3.3]{Su19}.

Such an approach is, however, inferior to the one with Tate twist in the case of
even~$i$,
at least as far as the component groups are considered. For example, the algebraic monodromy group might be
$[\SO_3(\bbQ_l)]^2$.
Then,~working without Tate twist, one would find, at first,
$\G[\SO_3(\bbQ_l)]^2$.
In~the next step, however, this leads to
$\G[\SO_3(\bbQ_l)]^2 \cap \O_6(\bbQ_l) = [\SO_3(\bbQ_l)]^2 \cup [\O_3^-(\bbQ_l)]^2$,
in which, all of a sudden, a second component~appears. In~other words, some of the information has been~lost.
\end{rem}

\subsubsection*{The Sato--Tate conjecture}

The set
$\Cl(\ST^i(X))$
of the conjugacy classes of elements
of~$\ST^i(X)$
naturally carries the quotient topology with respect to the canonical map
$\pi\colon \ST^i(X) \to \Cl(\ST^i(X))$.
As~$\smash{\tr|_{\ST^i(X)}\colon \ST^i(X) \to \bbC}$
is a continuous class function, it induces a continuous map
$\tr'\colon \Cl(\ST^i(X)) \to \bbC$
satisfying
$\smash{\tr' \!\circ\! \pi = \tr|_{\ST^i(X)}}$.
One~equips
$\Cl(\ST^i(X))$
with the measure
$\pi_* \mu_\Haar$,
for
$\mu_\Haar$
the normalised Haar measure
on~$\ST^i(X)$.

Moreover,~for an arbitrary
$p \in \bbP \setminus S$,
one~puts
$$\xi_p := \iota(\varrho^i_{X,l}(\Frob_p)) \in G_{X,l}^{i,\Zar}(\bbC) \,.$$
This~element is uniquely determined, up to conjugation.
Write~$\xi_p^\sss$
for the semisimple part
of~$\xi_p$,
according to the Jordan decomposition~\cite[Theorem~I.4.4]{Bo}. Then~all eigenvalues
of~$\smash{\xi_p^\sss}$
are of absolute
value~$1$.
Thus, the group
$\smash{\langle \xi_p^\sss\rangle \subset G_{X,l}^{i,\Zar}(\bbC)}$
has a compact~closure. By~\cite[Theorem IV.3.5]{OV},
$\smash{\langle \xi_p^\sss\rangle}$
is, up to conjugation, contained
in~$\ST^i(X)$.
Let,~finally,
$$x_p \in \Cl(\ST^i(X))$$
be the conjugacy class
of~$\xi_p^\sss$.
Lemma~\ref{cl_inj} below shows that
$x_p$~is
well-defined.

\begin{conj}[The Sato--Tate conjecture]
Let\/~$X$
and\/~$i$
be as~above. Then~the sequence\/
$(x_p)_{p \in \bbP \setminus S}$,
for\/~$p$
running through the good primes in their usual order, is equidistributed with respect
to\/~$\pi_* \mu_\Haar$.
In~other words, the sequence\/
$\big( \frac1{\#\{q \in \bbP \setminus S \mid q \le p\}} \!\! \sum\limits_{\genfrac{}{}{0pt}{}{q \in \bbP \setminus S}{q \le p}} \!\!\delta_{x_q} \big)_{p\in\bbP}$
of measures on\/
$\Cl(\ST^i(X))$
converges weakly
versus\/~$\pi_* \mu_\Haar$.\vspace{-2mm}
\end{conj}

\begin{rems}
\label{rems_allg}
\begin{iii}
\item
(Equidistribution on the component group.)
In particular, the Sato--\-Tate conjecture claims equidistribution among the components of
$\ST^i(X)$.

More~precisely, let
$\mu_\un$
be the uniform probability measure on the component~group
$\ST^i(X)/(\ST^i(X))^0$.
Moreover, let
$\pi_c\colon \ST^i(X)/(\ST^i(X))^0 \to \Cl(\ST^i(X)/(\ST^i(X))^0)$
be the canonical map and
$\kappa_c\colon \Cl(\ST^i(X)) \to \Cl(\ST^i(X)/(\ST^i(X))^0)$
the map between conjugacy classes induced by the projection
$\kappa\colon \ST^i(X) \to \ST^i(X)/(\ST^i(X))^0$.
Then
$(\kappa_c(x_p))_{p \in \bbP \setminus S}$
is asserted to be equidistributed with respect
to\/~$(\pi_c)_* \mu_\un$.
Indeed,
$$(\kappa_c)_* \pi_* \mu_\Haar = (\kappa_c \!\circ\! \pi)_* \mu_\Haar = (\pi_c \!\circ\! \kappa)_* \mu_\Haar = (\pi_c)_* \kappa_* \mu_\Haar = (\pi_c)_* \mu_\un \,. \vspace{2mm}$$

This part of the Sato--Tate conjecture is known to be true and can be shown as~follows. The image of
$\smash{\varrho^i_{X,l}\colon \Gal(\overline\bbQ/\bbQ) \to G_{X,l}^{i,\Zar}}$
is Zariski dense, hence the induced ho\-mo\-mor\-phism
$$\Gal(\overline\bbQ/\bbQ) \longrightarrow G_{X,l}^{i,\Zar}/(G_{X,l}^{i,\Zar})^0 \cong \ST^i(X)/(\ST^i(X))^0$$
is~surjective. The~kernel
$U \subseteq \Gal(\overline\bbQ/\bbQ)$
is an open subgroup, so corresponding under the Galois correspondence there is a finite extension field
$\smash{L_0 \supseteq \bbQ}$.
I.e.,
$\smash{\varrho^i_{X,l}}$
yields an isomorphism
$\smash{\Gal(L_0/\bbQ) \cong \ST^i(X)/(\ST^i(X))^0}$.
Consequently, the Chebotarev density theorem implies exactly what was~claimed.
\item
(The
$0$-dimensional
case.)
In~particular, the Sato--Tate conjecture is trivially true when
$(\ST^i(X))^0$
is the trivial~group. For~example, this holds for
$X$
of
dimension~$0$
and~$i=0$.
Indeed, then
$\smash{H^0_\et(X_{\overline\bbQ}, \bbQ_l) \cong \bbQ_l^{\#\pi_0(X_{\overline\bbQ})}}$
and
$\Gal(\overline\bbQ/\bbQ)$
acts simply by~permuting the direct~summands. Consequently,~the algebraic monodromy group
$\smash{G_{X,l}^{0,\Zar}}$
must be~finite.
\item
(Modularity.)
For a representation
$\varrho\colon \ST^i(X) \to \GL_d(\bbC)$,
consider the Artin type
$L$-function
$$L(\varrho,s) := \prod_{p \in \bbP \setminus S} \frac1{\det(1-\varrho(x_p) p^{-s})} \,,$$
which is clearly holomorphic
for~$\Re s > 1$.
Assume~that, for every irreducible, con\-tin\-u\-ous representation
$\varrho \neq 1$
of~$\ST^i(X)$,
the function
$L(\varrho,s)$
extends to the closed half plane
$\Re s \geq 1$
as a continuous function not having any zeroes (or poles). Then~the Sato--Tate conjecture is known to hold
for~$X$
and~$i$
\cite[\S I.A.2, Theorem~2]{Se68}.
\item
(Cohomology.)
The group
$\Gal(\overline\bbQ/\bbQ)$
is compact and hence carries a normalised Haar measure~itself. The~conjugacy classes of
$\Frob_p$,
for~$p$
running through the primes in their usual order, are equidistributed with respect to this Haar measure, according to the Chebotarev density~theorem. As~the representation
$\smash{\varrho^i_{X,l}}$
is continuous, the image
$\smash{\im \varrho^i_{X,l}}$
is compact, and the conjugacy classes of
$\smash{\varrho^i_{X,l}(\Frob_p)}$
are equidistributed with respect to the normalised Haar measure on that
\mbox{$l$-adic}
Lie~group. This
\mbox{$l$-adic}
kind of equidistribution is certainly of interest. For example, it was studied in detail, for elliptic curves, by J.-P.\ Serre in~\cite{Se72}. However,~as the embedding
$\bbQ_l \hookrightarrow \bbC$
chosen is discontinuous, it does not seem to have any implications towards the Sato--Tate conjecture. A~cohomology theory with coefficients
in~$\bbC$
that provides a continuous
\mbox{$\Gal(\overline\bbQ/\bbQ)$-action}
would certainly help. But, of course, we have nothing of this kind at our~disposal.
\end{iii}
\end{rems}

\begin{rem}
When~$X$
is a
$K3$~surface,
which is the situation we are interested in in this article,
$\xi_p$
is known to be semisimple, for every good prime
$p \neq l$~\cite[Corollaire 1.10]{De81}.
The step of taking the semisimple part is then superfluous.
\end{rem}

The~Sato--Tate conjecture immediately yields the following prediction for the distribution of the Frobenius~traces.

\begin{conj}[The Frobenius trace distribution]
Let\/~$X$
and\/~$i$
be as~above. Then~the sequence\/
$\smash{\Tr(\Frob_p)_{p \in \bbP \setminus S}}$,
for\/~$p$
running through the good primes in their usual order, is equidistributed with respect
to\/~$(\tr|_{\ST^i(X)})_* \mu_\Haar$.
In~other words, the sequence\/
$$\textstyle \big(\frac1{\#\{q \in \bbP \setminus S \mid q \le p\}} \!\!\sum\limits_{\genfrac{}{}{0pt}{}{q \in \bbP \setminus S}{q \le p}}\!\! \delta_{\Tr(\Frob_q)}\big)_{p\in\bbP}$$
of measures
on\/~$\bbR$
is convergent in the weak sense
versus\/~$(\tr|_{\ST^i(X)})_* \mu_\Haar$.\medskip

\noindent
{\bf Proof}
{\em (assuming the Sato--Tate conjecture)}%
{\bf .}
{\em
Taking the image measure under a continuous map commutes with weak convergence, cf.\ \cite[section 13.4, probl\`eme~8]{Di}. Hence,~the Sato--Tate conjecture implies~that
\begin{align*}
\textstyle \frac1{\#\{q \in \bbP \setminus S \mid q \le p\}} \!\!\sum\limits_{\genfrac{}{}{0pt}{}{q \in \bbP \setminus S}{q \le p}}\!\! \delta_{\tr'(x_q)} = \tr'_* \big(\frac1{\#\{q \in \bbP \setminus S \mid q \le p\}} \!\!\sum\limits_{\genfrac{}{}{0pt}{}{q \in \bbP \setminus S}{q \le p}}\!\! \delta_{x_q}\big) &\to \tr'_*(\pi_* \mu_\Haar) \\[-5.5mm]
&= (\tr' \!\circ\! \pi)_* \mu_\Haar = (\tr|_{\ST^i(X)})_* \mu_\Haar \,.
\end{align*}
But
$\smash{\tr'(x_q) = \tr'(\pi(\xi_q)) = \tr(\xi_q) = \tr(\iota(\varrho^i_{X,l}(\Frob_p)))}$,
for every prime
number~$q$.
And
$\smash{\tr(\iota(\varrho^i_{X,l}(\Frob_p))) = \Tr(\varrho^i_{X,l}(\Frob_p))}$,
which is usually denoted shortly
as~$\Tr(\Frob_p)$.
}
\eop
\end{conj}

\subsubsection*{A related conjecture -- Lang--Trotter for general varieties}

Let,~as before,
$X$
be a smooth, projective variety
over~$\bbQ$,
and
$S$
a finite set of primes, outside of which
$X$
has good~reduction. Fix some
$i\in\bbN$.
Then,~for any
$a\in\bbZ$,
one may ask for the asymptotics~of
\begin{equation}
\label{LT}
N^i_{X,a}(p) := \#\{q\in\bbP \setminus S \mid q \leq p, \Tr(\Frob_q\colon H_\et^i(X_{\overline\bbQ}, \bbQ_l) \righttoleftarrow) = a \} \,,
\end{equation}
for~$p\to\infty$.
Note~that, as the cohomology vector space without Tate twist is considered, the traces of the Frobenii in~(\ref{LT}) are automatically~integers.

Formula~(\ref{tate}) shows that, in the notation used above, the second condition in the definition
of~$N^i_{X,a}(p)$
means
$\smash{\Tr(\Frob_q) = \frac{a}{q^{i/2}}}$,
or
$$\textstyle \smash{\Tr(\Frob_q) \in \left[ \frac{a-\frac12}{q^{i/2}}, \frac{a+\frac12}{q^{i/2}} \right]} \,.$$
Thus,~assuming that the Sato--Tate measure
$(\tr|_{\ST^i(X)})_* \mu_\Haar$
has a density function whose limit
for~$t\to0$
exists and is positive, it seems reasonable to expect, at least for
$a \neq 0$,
that there is a constant
$C^i_{X,a} \geq 0$
such~that
$$N^i_{X,a}(p) \sim C^i_{X,a} \!\cdot\! \!\sum_{\genfrac{}{}{0pt}{}{q \in \bbP \setminus S}{q \le p}} \!\frac1{q^{i/2}} \,.$$
I.e., that
$N^i_{X,a}(p) = O(1)$,
for~$i>2$,
and
\begin{equation}
\label{LTconj}
N^i_{X,a}(p) \sim
\left\{
\begin{array}{ll}
C^i_{X,a} \!\cdot\! \frac{2\sqrt{p}}{\log p} \,, & {\rm ~for\;} i=1 \,, \\
C^i_{X,a} \!\cdot\! \log\log p \,,               & {\rm ~for\;} i=2 \,. \\
\end{array}
\right.
\end{equation}
In~the case that the density is continuous and non-vanishing
at~$0$,
one might hope for the same
when~$a=0$.
This~was formulated first, as a conjecture for elliptic curves
and~$i=1$,
by S.~Lang and H.~Trotter \cite{LT}.

\subsubsection*{The case of a
$K3$~surface}

For~$K3$~surfaces
and~$i=2$,
which is the case we are interested in in this article, it seems that such a conjecture has not been explicitly stated~before. Note,~however, the somewhat optimistic discussion on page~2 of the article~\cite{CT} of E.~Costa and~Yu.\ Tschinkel.

Anyway,~we are very reluctant to claim evidence in any nontrivial situation, as the experiments described in this article involve a search bound
of~$10^8$,
which is a bit too low in order to detect double logarithmic~growth.

There~are, however, cases, in which (\ref{LTconj}) is trivially~true. For~instance, let
$X$
be one of the surfaces from Examples~\ref{ex1} to~\ref{ex4},~below. Then,~for
$a$~even,
(\ref{LTconj})~holds
with~$C^2_{X,a} = 0$.
Indeed, a double cover
of~$\smash{\Pb^2_{\bbF_{\!p}}}$,
ramified over six
\mbox{$\bbF_{\!p}$-rational}
lines in general position, always has an odd number of
\mbox{$\bbF_{\!p}$-rational}
points, simply because the branch locus~has. For~a refinement of this argument taking
$(a \bmod 4)$
into consideration, cf.~\cite{EJ22}.

\begin{rem}
One might expect an asymptotics similar to (\ref{LTconj}) for the primes of reduction to geometric Picard
rank~$22$.
Again, this is certainly subject to restrictions, such as the occurrence of a continuous density function for the Sato--Tate measure. Note,~for instance, that the
$K3$~surfaces
presented in \cite[Examples 2.6.5 and~2.6.7]{CEJ} reduce to geometric Picard
rank~$22$
at exactly half the~primes.
\end{rem}

\subsubsection*{A result from the theory of Lie~groups}

Due to the lack of a suitable reference, we include the following purely Lie-the\-o\-retic lemma.

\begin{lem}
\label{cl_inj}
Let\/~$G$
be a faithfully representable complex Lie group and\/
$K \subset G$
a maximal compact subgroup. Then the natural homomorphism\/
$\Cl K \to \Cl G$
between conjugacy classes of elements is injective.\medskip

\noindent
{\bf Proof.}
{\em
Let
$k_1,k_2 \in K$
be two elements that are conjugate as elements
of~$G$.
We~have to show that
$k_1$
and~$k_2$
are conjugate
in~$K$.

According to the decomposition theorem \cite[Theorem~4.43]{Le},
$G \cong G' \!\ltimes\! R$
is isomorphic to a semidirect product of two closed subgroups,
$R$
being simply connected and solvable and
$G'$
being reductive. A~simply connected solvable group has the trivial group as its maximal compact subgroup~\cite[Corollary~1.126]{Kn}. Thus, the quotient homomorphism
$\pi\colon G \twoheadrightarrow G/R$
maps
$K$
isomorphically onto the maximal compact subgroup
of~$G/R$.
Moreover, the elements
$\pi(k_1)$
and~$\pi(k_2)$
are clearly conjugate
in~$G/R \cong G'$.

It therefore suffices to assume
$G$
as being reductive. Then
$G$
coincides with the complexification
of~$K$~\cite[Theorem~4.31]{Le}
and there is the Cartan decomposition
$G = K \!\cdot \exp(i \Lie K)$,
cf.~\cite[Theorem~6.31.c)]{Kn}. By~assumption, there exist some
$k \in K$
and
$X \in \Lie K$
such that
$\exp(-iX)k^{-1} \!\cdot\! k_1 \!\cdot\! k\exp(iX) = k_2 \in K$.
As~$K$
is fixed under the Cartan involution
$\Theta$,
this~yields
$$\exp(iX) \underline{k} \exp(-iX) = \Theta(\exp(-iX) \underline{k} \exp(iX)) = \exp(-iX) \underline{k} \exp(iX) \,,$$
for
$\underline{k} := k^{-1}k_1k$.
I.e.,
$\exp(2iX) \in \C_G(\underline{k})$.
Consequently,
$\exp(2niX) \in \C_G(\underline{k})$
for every
$n \in \bbZ$,
which implies
$\exp(tiX) \in \C_G(\underline{k})$
for every
$t \in \bbR$
\cite[Lemma~1.142]{Kn}. In~particular,
$\exp(iX) \in \C_G(\underline{k})$
showing
$k^{-1} k_1 k = k_2$,
as~required.%
}%
\eop
\end{lem}

\section{Trace distributions for compact Lie groups}
\label{Lie}

\subsubsection*{Moment sequences}

Given a connected compact Lie group
$S \subset \GL_d(\bbC)$,
the trace
$\tr|_S\colon S \to \bbC$
is a continuous class~function. Let~us assume that
$\tr|_S$
is real-valued. Then,~for the
\mbox{$n$-th}
moment
$$\E_S[\tr^n] = \int_\bbR \!x^n\, d(\tr|_S)_*\mu_\Haar(x) = \int_S \tr^n(s)\, d\mu_\Haar(s) \,,$$
for
$n\in\bbN$,
one has the Weyl integration formula \cite[Theorem 8.60]{Kn},
$$\E_S[\tr^n] = \frac1{\#W} \int_T \tr^n(t) \det(\id-\Ad_{S/T}(t^{-1})) \,dt = \frac1{\#W} \int_T \tr^n(t)\! \prod_{\alpha\in\Delta^+}\!\! |1-t^\alpha|^2 \,dt \,.$$
Here,
$T$
denotes the maximal torus
of~$S$,
$W$
the Weyl group,
$\Delta^+$
a system of positive roots, and
$\Ad\colon T \to \GL(S/T)$
the adjoint representation. The integrand is a trigonometric polynomial, so, for
each~$n$,
the integral may be computed~exactly.

The~particular compact Lie groups mentioned in Table~1 are related to the examples of
$K3$~surfaces
that we present in this article, cf.\ Section~\ref{sec_exp}.\medskip\smallskip

\subsubsection*{Examples~A}
~\vspace{-2mm}

\begin{table}[H]
\begin{center}
\begin{tabular}{|l||c|c|c|c|c|}
\hline
$S$                      & Root        & $\dim(T)$ & $\dim(S/T)$ & Moment sequence          & Label in    \\
                         & system      &           &             &                          & \cite{OEIS} \\\hline\hline
$\SO_2(\bbR) \cong \U_1$ & $\emptyset$ & 1         &  0          & 1, 0, 2, 0, 6, 0, 20,    & A126869     \\
                         &             &           &             & 0, 70, 0, 252            &             \\\hline
$\SO_3(\bbR)$            & $A_1$       & 1         &  2          & 1, 0, 1, 1, 3, 6, 15,    & A005043     \\
                         &             &           &             & 36, 91, 232, 603         &             \\\hline
$\SO_5(\bbR)$            & $B_2$       & 2         &  8          & 1, 0, 1, 0, 3, 1, 15,    & A095922     \\
                         &             &           &             & 15, 105, 190, 945        &             \\\hline
$\SO_6(\bbR)$            & $A_3$       & 3         & 12          & 1, 0, 1, 0, 3, 0, 16,    & A247591     \\
                         &             &           &             & 0, 126, 0, 1296          &             \\\hline
$\U_3$                   & $A_2'$      & 3         &  6          & 1, 0, 2, 0, 12, 0, 120,  & A245067     \\
                         &             &           &             & 0, 1610, 0, 25956        &             \\\hline
\end{tabular}
\end{center}\vspace{-2mm}

\caption{Moments of the trace distributions for some connected compact Lie groups}\vspace{-2mm}

\end{table}

\begin{rems}
\begin{iii}
\item
$\U_3$
is reductive, but not semisimple. Thus, the root system
$A_2$
occurs together with a
\mbox{$3$-dimensional}
torus.
\item
In the
\mbox{$\SO_N(\bbR)$-cases},
we consider the naive trace function,
$\tr\colon A \mapsto \Tr(A)$.
In~the
\mbox{$\U_N$-cases},
however, we put
$\smash{\tr(A) := \Tr \genfrac(){0pt}{}{A \,0}{\,0 \,\overline{A}} = 2 \Re \Tr(A)}$.
This is compatible with the embedding
$\GL_N(\bbC) \hookrightarrow \SO_{2N}(\bbC)$,
$\smash{A \mapsto \genfrac(){0pt}{}{\,E\;\;\;\;\;E}{iE \,-iE}^{-1} \genfrac(){0pt}{}{A \;\;\;\;0\;\;\;\;\;\;}{0 \,(A^t)^{-1}} \genfrac(){0pt}{}{\,E\;\;\;\;\;E}{iE \,-iE}}$,
of complex Lie~groups, which is relevant here. Cf.\ Theorem~\ref{neutr_comp}.ii), below.
\item
The moment sequences for the classical groups, such as those presented in Table~1, have beautiful combinatorial interpretations~\cite[Theorem~3.16]{Me}.
\end{iii}
\end{rems}

\subsubsection*{Examples~B}

For~$S = S_1 \times S_2$
and
$\tr(s_1,s_2) = \tr_{S_1}(s_1) + \tr_{S_2}(s_2)$,
according to Fubini, one~has
$$\E_S[\tr^n] = \sum_{i=0}^n \genfrac(){0pt}{}{n}{i} \E_{S_1}[\tr_{S_1}^i] \E_{S_2}[\tr_{S_2}^{n-i}] \,.$$
For~two Lie groups of this kind, the moment sequences are of interest for~us.

\begin{table}[H]
\begin{center}
\begin{tabular}{|l||c|c|c|c|}
\hline
$S$                              & Root        & $\dim(T)$ & $\dim(S/T)$ & Moment sequence          \\
                                 & system      &           &             &                          \\\hline\hline
$[\SO_3(\bbR)]^2$                & $[A_1]^2$   & 2         &  4          & 1, 0, 2, 2, 12, 32, 140, \\
                                 &             &           &             & 534, 2324, 10112, 46008  \\\hline
$[\SO_2(\bbR)]^3 \cong [\U_1]^3$ & $\emptyset$ & 3         &  0          & 1, 0, 6, 0, 90, 0, 1860, \\
                                 &             &           &             & 0, 44730, 0, 1172556     \\\hline
\end{tabular}
\end{center}\vspace{-2mm}

\caption{Moments of the trace distributions in direct product cases}
\end{table}

\subsubsection*{Example~C}

There is a version of the Weyl integration formula for disconnected compact Lie groups \cite[Proposition~2.3]{We}, of which we made use in the case
of~$\O_6(\bbR)$.

\begin{table}[H]
\begin{center}
\begin{tabular}{|c||c|c|c|c|}
\hline
Component                            & $\dim(T^\prime)$ & $\dim(S/T^\prime)$ & Moment sequence        & Label in \cite{OEIS}  \\\hline\hline
$\O_6^-(\bbR) = {}\hspace{.93cm}$      & 2                & 13                 & 1, 0, 1, 0, 3, 0, 14,  & A138349               \\
$\O_6(\bbR) \setminus\! \SO_6(\bbR)$ &                  &                    & 0, 84, 0, 594, 0, 4719 &                       \\\hline
\end{tabular}
\end{center}\vspace{-2mm}

\caption{Moments of the trace distribution for a non-neutral component}
\end{table}

\subsubsection*{Plotting the density of the trace distribution}

We~used two rather different approaches to plot the density of a trace distribution.

\begin{iii}
\item
(The moments based approach.)
We~split the support interval of the density into subintervals of equal lengths. Then we compute a cubic spline function, using the subdivision chosen, that has the same moments as the distribution. We~work with the first 20 to 35~moments.
\item
(The numerical integration approach.)
For every continuous function
$g\colon \bbR \to \bbR$,
one has, again by the Weyl integration formula,
\begin{align*}
(\tr|_S)_*\mu_\Haar(g) = \int_\bbR g(x)\, d(\tr|_S)_*\mu_\Haar(x) = \int_S g&(\tr(s))\, d\mu_\Haar(s) = \\[-2mm]
&= \frac1{\#W} \int_T g(\tr(t))\! \prod_{\alpha\in\Delta^+}\!\! |1-t^\alpha|^2 \,dt \,.
\end{align*}
Thus,~the density function
$\den_S$
of the distribution
$(\tr|_S)_*\mu_\Haar$
can be described as follows.

For~$x\in\bbR$,
put
$T_x := \{t\in T \mid \tr(t) = x\}$.
At~least for
$x \in \bbR \setminus \bbZ$
in the range
of~$\tr$,
this is a submanifold
of~$T$
of
codimension~$1$.
Then
$$\den_S(x) = \frac1{\#W} \int_{T_x} \frac{\prod_{\alpha\in\Delta^+}\! |1-t^\alpha|^2}{|\frac{\partial\,\tr}{\partial\vec{\frakn}}(t)|} \,dt \,,$$
for
$\vec{\frakn}$
a normal vector of
length~$1$,
as~usual.
Using numerical integration, we evaluate these
\mbox{$(\dim(T)-1)$-dimensional}
integrals at a sufficiently high~precision.
\end{iii}

\begin{rem}
Plotting the densities with either approach, the results visually coincide. The only exception is the case
of~$[U_1]^3$.
In~fact, the moments based approach does not work properly
for~$[U_1]^3$.
As~the density in this case is not a
$C^2$~function,
an approximation by cubic splines is not appropriate. The plot shows oscillations that increase when refining the subdivision. In~this case, the more naive numerical integration approach has to be~applied.
\end{rem}

\subsubsection*{Explicit formulas for some of the density functions}

Instead of either of the two approaches, one might ask for explicit formulas for the density functions, analogous to those in the elliptic curves~case. Unfortunately, the possibilities for such an approach are very limited.

\begin{iii}
\item
For~$\SO_2(\bbR) \cong \U_1$,
the density function is
$\smash{\frac1{\pi\sqrt{4-t^2}}}$
on~$[-2,2]$.
This~is well-known, since it concerns the case of CM elliptic~curves.

For~$\SO_3(\bbR)$,
the density function is given
on~$[-1,3]$
by~$\smash{\frac1{2\pi}\sqrt{\frac{3-t}{1+t}}}$,
as is shown by an elementary calculation (cf.\ \cite[Formula~(4.85)]{KM}).

Let~us note that
for~$\SO_4(\bbR)$,
which we do not need any further in this article, an explicit formula for the density function is known, too~\cite[Theorem~8.7]{EP}.
For~$\SO_5(\bbR)$,
G.\ Lachaud \cite[Proposition~8.5 and Remark~8.6]{La} showed that the density function~is
$$-\frac{1}{24\pi^2} \textstyle \left((t-1)(3t^2+2t+43) \K\!\big(1-\frac{(t-1)^2}{16}\big) - 4(t^2+22t-7) \E\!\big(1-\frac{(t-1)^2}{16}\big)\right)$$
on~$[-3,5]$.
Here,~$\K$
and~$\E$
denote the complete elliptic integrals of the first and second kinds, given by
$\smash{\K(t) = \int_0^{\pi/2} \!d\theta / \sqrt{1-t \sin^2\!\theta}}$
and
$\smash{\E(t) = \int_0^{\pi/2} \!\sqrt{1-t \sin^2\!\theta} \,d\theta}$,
respectively. Both~are holomorphic functions on
$\bbC \,\setminus\, [1,\infty)$
and positively real-valued
on~$(-\infty,1)$.

We~do not know, however, of an explicit formula for the density function in
the case
of~$\SO_6(\bbR)$.
Neither do we
for~$\U_3$.
Our~attempts to replace the numerical integration by symbolic methods turned out unsuccessful for these Lie~groups. At~least
for~$\SO_6(\bbR)$,
one should probably expect a more complicated answer than
for~$\SO_5(\bbR)$.
\item
For~$[\SO_3(\bbR)]^2$,
the convolution of the density function
for~$\SO_3(\bbR)$
with itself is asked~for. A~calculation in {\tt maple} yields
$$\frac1{4\pi^2} \textstyle \left((2-t) \K\!\big(1-\frac{(t-2)^2}{16}\big) + 4 \E\!\big(1-\frac{(t-2)^2}{16}\big)\right)$$
on~$[-2,6]$.

Moreover,
for~$[\U_1]^2$,
the density
on~$[-4,4]$
turns out to be given~by
\begin{equation}
\label{U_1^2}
\frac1{2\pi^2} \textstyle \K(1-\frac{t^2}{16}) \,.
\end{equation}

Thus,~the density function
for~$[\U_1]^3$
is the convolution of~(\ref{U_1^2}) with
$\smash{\frac1{\pi\sqrt{4-t^2}}}$
on~$[-2,2]$.
An~explicit formula is known,~too. Indeed,~from \cite[Formulas (7.10) to~(7.14)]{Jo}, one finds
\begin{equation}
\label{U_1^3}
- \frac2{\pi^3 |t|} \Ima \!\big(\sqrt{4-v(t)}\, \sqrt{1-v(t)}\,\K(m_+(t)) \,\K(m_-(t))\big)
\end{equation}
on~$[-6,6]$,
for
$v(t) := \frac{20 - t^2 + \sqrt{(4-t^2)(36-t^2)}}8$
and
$m_{\pm}(t) := \frac{2 \pm v(t)\sqrt{4-v(t)} - (2-v(t))\sqrt{1-v(t)}}4$.
Here,
$(4-t^2)(36-t^2)$
is real, so one takes the natural branch of the square root on the upper half~plane. On~the other hand, for the square roots of
$(4-v(t))$
and
$(1-v(t))$,
the natural branch of the square root on the lower half~plane is~taken. Moreover,
at~$t = \pm2$,
it happens that
$m_+(t)$
crosses the branch cut
of~$\K$,
but the imaginary part in~(\ref{U_1^3}) is nevertheless~continuous.
\item
Finally,~for
$\O_6^-(\bbR)$,
the density function may be written
on~$[-4,4]$~as
$$-\frac1{60\pi^2} \textstyle \left((8t^4+128t^2) \K(1-\frac{t^2}{16}) - (t^4+224t^2+256) \E(1-\frac{t^2}{16})\right) \,.$$
This~is essentially the formula given in \cite[Corollary~5.6]{La}. Observe~that the author considers an equivalent distribution. It~is interesting to note that he provides many more expressions for the same density in terms of other special~functions.
\end{iii}

\section{$K3$~surfaces}

\subsubsection*{Decomposition of cohomology}

Let~$X$
be a
$K3$~surface
over~$\bbQ$.
One then calls
$H_\tr := (H_\alg)^\perp \subset H^2(X(\bbC), \bbQ)$
the transcendental part of the cohomology. Here,
$$H_\alg := \im(c_1\colon \Pic X(\bbC) \to H^2(X(\bbC), \bbQ)) \,.$$
As~a quadratic space,
$H_\alg$~is
non-degenerate. Indeed,~if
$\calL \in \Pic X(\bbC)$,
$\calL \not\cong \calO_{X(\bbC)}$,
had intersection
number~$0$
with every element
of~$\Pic X(\bbC)$
then either
$\calL$
or~$\calL^\vee$
would have a non-trivial section~\cite[Proposition~VIII.3.7.i)]{BHPV}. And~hence
$\calL \!\cdot\! [H] \neq 0$,
for~$H$
the hyperplane section, a~contradiction. Consequently,
$H_\tr = (H_\alg)^\perp$
is non-degenerate,~too.

\begin{nota}
Write
$r := \rk\Pic X(\bbC)$.
Then~$\dim_\bbQ H_\tr = 22-r$.
\end{nota}

In
\mbox{$l$-adic}
cohomology, one puts
$\smash{H_{l,\alg} := \im(c_1\colon \Pic X_{\overline\bbQ} \!\otimes_\bbZ\! \bbQ_l \to H^2_\et(X_{\overline\bbQ}, \bbQ_l(1)))}$
and
$H_{l,\tr} := (H_{l,\alg})^\perp$.
Then, under the standard comparison isomorphism~\cite[Expos\'e~XVI, Th\'eor\`eme 4.1]{SGA4},
$H_{l,\tr} \cong H_\tr \!\otimes_\bbQ\! \bbQ_l(1)$
and
$H_{l,\alg} \cong H_\alg \!\otimes_\bbQ\! \bbQ_l(1)$.
The~representation
$\varrho^2_{X,l}$
maps
$H_{l,\alg}$~to
itself and therefore
$H_{l,\tr}$~to
itself,~too. Thus,
$\varrho^2_{X,l}$
splits into the direct sum of the two sub-representations
$$\varrho^2_{X,l,\alg}\colon \Gal(\overline\bbQ/\bbQ) \to \O(H_{l,\alg}) \quad{\rm and}\quad \varrho^2_{X,l,\tr}\colon \Gal(\overline\bbQ/\bbQ) \to \O(H_{l,\tr}) \,.$$
The~image
of~$\varrho^2_{X,l,\alg}$
is a finite group
$G_{\Pic} \cong \Gal(K_{\Pic}/\bbQ)$,
for
$K_{\Pic}$
the splitting field
of~$\Pic X_{\overline\bbQ}$.

\begin{defi}
\begin{iii}
\item
We~call
$\smash{G_{X,l,\tr}^{2,\Zar} := \overline{\im(\varrho^2_{X,l,\tr})} \subseteq \O(T_l)}$
the {\em transcendental part of the algebraic monodromy group\/}
of~$X$.
\item
Moreover,~the maximal compact subgroup
of~$\smash{G_{X,l,\tr}^{2,\Zar}(\bbC)}$
is called the {\em transcendental part of the Sato--Tate group\/}
of~$X$
and denoted
by~$\ST^2_\tr(X)$.
\end{iii}
\end{defi}

\begin{lem}
\begin{abc}
\item
For the neutral components, one has
$$(G_{X,l,\tr}^{2,\Zar})^0 = (G_{X,l}^{2,\Zar})^0 \,,
\;\;
(G_{X,l,\tr}^{2,\Zar}(\bbC))^0 = (G_{X,l}^{2,\Zar}(\bbC))^0 \,
\;\;{\it and}\;\;
(\ST^2_\tr(X))^0 = (\ST^2(X))^0 \,.$$
\item
Concerning the component groups,
$$G_{X,l,\tr}^{2,\Zar}/(G_{X,l}^{2,\Zar})^0 = G_{X,l,\tr}^{2,\Zar}(\bbC)/(G_{X,l}^{2,\Zar}(\bbC))^0 \cong \ST^2_\tr(X)/(\ST^2(X))^0 \,.$$
\end{abc}\smallskip

\noindent
{\bf Proof.}
{\em
a)
The first equality is a direct consequence of the fact that
$
G_{\rm Pic}$
is~finite. The second one follows immediately from the first, and, finally, the third is obtained taking the maximal compact subgroup on either~side.\smallskip

\noindent
b)
follows from the standard facts that the maximal compact subgroup of a Lie group meets every connected component, and that the maximal compact subgroup of a connected Lie group is~connected.
}
\eop
\end{lem}

\begin{lem}
The homomorphism
$$G_{X,l}^{2,\Zar}/(G_{X,l}^{2,\Zar})^0 \longrightarrow G_{\rm Pic} \oplus G_{X,l,\tr}^{2,\Zar}/(G_{X,l}^{2,\Zar})^0$$
induced by the decomposition, is a subdirect product. I.e., it is an injection, but the projections to either summand are~surjective.\medskip

\noindent
{\bf Proof.}
{\em
One has that
$\smash{G_{X,l}^{2,\Zar}}$
is the image of
$\smash{\varrho^2_{X,l} = \varrho^2_{X,l,\alg} \oplus \varrho^2_{X,l,\tr}}$,
while
$G_{\rm Pic}$
and
$\smash{G_{X,l,\tr}^{2,\Zar}}$
are the images of the direct~summands. Therefore, the decomposition induces a homomorphism
$\smash{G_{X,l}^{2,\Zar} \to G_{\rm Pic} \oplus G_{X,l,\tr}^{2,\Zar}}$
that is a subdirect product. The assertion follows immediately from~this.
}
\eop
\end{lem}

In certain cases, the trace of the Frobenius
on~$\ST^2_\tr(X)$
is related to the point count on a singular model of the
$K3$~surface~$X$.

\begin{lem}
\label{bl_down}
Let\/~$X$
be a\/
$K3$~surface
over\/~$\bbQ$
and\/
$\pr\colon X \to X'$
a birational~morphism.
Write\/~$\smash{r := \rk\Pic X_{\overline\bbQ}}$
and let\/
$r_0$
be the number of\/
\mbox{$(-2)$-curves}
blown down under\/
$\smash{\pr_{\overline\bbQ}}$.
Suppose~that these generate\/
$\smash{\Pic X_{\overline\bbQ}}$,
together with\/
$(r-r_0)$
further linearly independent classes that are defined
over\/~$\bbQ$.\smallskip

\noindent
Then,~for every
prime\/~$p$
of good~reduction,
$$\#\calX'_p(\bbF_{\!p}) = 1 + p(r-r_0) + p \!\cdot\! \Tr(\varrho^2_{X,l,\tr}(\Frob_p)) + p^2 \,.$$
{\bf Proof.}
{\em
For a
$K3$~surface,
$\smash{H^1_\et(X_{\overline\bbQ}, \bbQ_l) = H^3_\et(X_{\overline\bbQ}, \bbQ_l) = 0}$,
while
$\smash{H^0_\et(X_{\overline\bbQ}, \bbQ_l)}$
and
$\smash{H^4_\et(X_{\overline\bbQ}, \bbQ_l)}$
are one-dimensional. Hence, in view of~(\ref{tate}), formula~(\ref{lefschetz}) specialises~to
\begin{align*}
\#\calX_p(\bbF_{\!p}) &= 1 + p \!\cdot\! \Tr(\varrho^2_{X,l}(\Frob_p)) + p^2 \\
 &= 1 + p \!\cdot\! \Tr(\varrho^2_{X,l,\alg}(\Frob_p)) + p \!\cdot\! \Tr(\varrho^2_{X,l,\tr}(\Frob_p)) + p^2 \,.
\end{align*}
Thus,~one has to show that
$\#\calX_p(\bbF_{\!p}) - \#\calX'_p(\bbF_{\!p}) + p(r-r_0) = p \!\cdot\! \Tr(\varrho^2_{X,l,\alg}(\Frob_p))$.

For~this, let us consider
the~$r_0$
points blown~up. These~are permuted
by~$\Frob_p$
and the difference
$\smash{\#\calX_p(\bbF_{\!p}) - \#\calX'_p(\bbF_{\!p})}$
may be written as
$p$~times
the number of fixed points. Which is the same as
$p \cdot\! \Tr(\varrho_\bl(\Frob_p))$,
for~$\varrho_\bl$
the corresponding permutation representation, a sub-representation
of~$\smash{\varrho^2_{X,l,\alg}}$.
As~the complement
of~$\varrho_\bl$
is, by assumption, trivial of
rank~$(r-r_0)$,
the claim~follows.
}
\eop
\end{lem}

\begin{ttt}[The neutral component of the Sato--Tate group]
\label{neut_comp}
The transcendental part of the cohomology is a pure
\mbox{weight-$2$}
Hodge structure
$H_\tr \subset H^2(X(\bbC), \bbQ)$.
Pure~Hodge structures of a fixed weight form an abelian category~\cite[Paragraphe~2.1.11]{De71}. The~endomorphisms
of~$H_\tr$,
as a Hodge structure, hence form a commutative
ring~$E$
with~$1$,
the {\em endomorphism ring\/}
of~$H_\tr$.
It~is well-known that
$E$
is, as long as
$K3$
surfaces are considered, always a~field.\medskip

\noindent
For~$X$
any
$K3$~surface
over~$\bbC$,
there are exactly three possibilities \cite[Theorem 1.6.a)]{Za}.

\begin{iii}
\item
One has
$E=\bbQ$.
This is the generic case.
\item
$E \supsetneqq \bbQ$
is a totally real field.
Then~$X$
is said to have {\em real multiplication (RM)}.

Put~$\delta := [E:\bbQ]$
and let
$e \in E$
be a primitive element. It~is known that
$e$~acts
on~$H_\tr$
as a self-adjoint linear map~\cite[Theorem 1.5.1]{Za}. Thus,
$H_\tr \!\otimes_\bbQ\! \bbC = H_1 \oplus\cdots\oplus H_\delta$
splits into eigenspaces that are mutually perpendicular. The~eigenspaces
$H_i$,
for
$i=1,\ldots,\delta$,
are not defined
over~$\bbQ$,
but form a single orbit under conjugation
by~$\Gal(\overline\bbQ/\bbQ)$.
In~particular,
$\smash{\dim_\bbQ H_1 = \cdots = \dim_\bbQ H_\delta = \frac{22-r}\delta}$.
Let~us note, in addition, that
each~$H_i$
is a simultaneous eigenspace for all elements
of~$E$.
\item
$E$~is
a CM field.
Then~$X$
is said to have {\em complex multiplication (CM)}.

Write
$E = E_0(\sqrt{-\tau})$,
for
$E_0 \subset E$
the maximal totally real subfield and
$\tau\in E_0$
a totally positive~element.
We~put
$\delta := [E_0:\bbQ]$.
Then, as above, the action
of~$E_0$
splits
$H_\tr \!\otimes_\bbQ\! \bbC = H_1 \oplus\cdots\oplus H_\delta$
into simultaneous eigenspaces, which are mutually perpendicular.

Under~the
action~$I$
of~$\sqrt{-\tau}\in E$,
each~$H_i = H_{i,+} \oplus H_{i,-}$,
for
$i = 1,\ldots,\delta$,
is split into two eigenspaces.
For~$v$
and~$w$
in the same eigenspace, one has~\cite[Theorem 1.5.1]{Za}
$$-\tau \!\cdot\! \langle v,w \rangle  = \langle Iv, Iw \rangle  = \langle v, -I^2w \rangle  = \tau \!\cdot\! \langle v,w \rangle \,,$$
which yields that the eigenspaces
$H_{i,+}$
and~$H_{i,-}$
are both~isotropic.
\end{iii}
\end{ttt}

\begin{rem}
The Hodge conjecture for
$(X \!\times\! X)(\bbC)$
implies that every endomorphism of
$H_\tr \subset H^2(X(\bbC), \bbQ)$
is induced by a correspondence
$S \subset (X \!\times\! X)(\bbC)$.
There~are two~issues.

\begin{iii}
\item
Such~a correspondence is clearly not~unique.
\item
As~a
$K3$~surface
does not carry a natural group structure, there is no reason to expect the endomorphisms
of~$H_\tr$
to be induced by self-morphisms
of~$X(\bbC)$.
\end{iii}
\end{rem}

\begin{theo}[Zarhin, Tankeev]
\label{neutr_comp}
Let\/~$X$
be a\/
$K3$~surface
over\/~$\bbQ$.
Moreover, let
$\smash{H_\tr \subset H^2(X(\bbC), \bbQ)}$
be the transcendental part of the cohomology,
and\/~$E$
its endomorphism~field.

\begin{iii}
\item
If\/~$E$
is totally real of
degree\/~$\delta$~then
$$\smash{(G_{X,l}^{2,\Zar}(\bbC))^0 \cong \SO(H_1) \times\cdots\times \SO(H_\delta) \cong [\SO_{\frac{22-r}\delta}(\bbC)]^\delta \,.}$$
For\/~$\delta=1$,
this includes the generic
case~$E=\bbQ$.
\item
If\/~$E$
is a CM~field of
degree\/~$2\delta$~then
$$\smash{(G_{X,l}^{2,\Zar}(\bbC))^0 \cong \O(H_1)_{(H_{1,+},H_{1,-})} \times\cdots\times \O(H_\delta)_{(H_{\delta,+},H_{\delta,-})} \cong [\GL_{\frac{22-r}{2\delta}}(\bbC)]^\delta \,.}$$
\end{iii}

\noindent
{\bf Proof.}
{\em
Due to the work of S.\,G.\ Tankeev~\cite{Tan90,Tan95}, together with \cite[Theorem 2.2.1]{Za}, one~has
$$(G_{X,l}^{2,\Zar}(\bbC))^0 \cong (\C_E(\O(H_\tr \!\otimes_\bbQ\! \bbC)))^0 \,.$$
Since a linear map commutes with the action
of~$E$
if and only if it maps each of the eigenspaces
$H_1, \ldots, H_\delta$,
or
$\smash{H_{1,+}, H_{1,-}, \ldots, H_{\delta,+}, H_{\delta,-}}$,
respectively, to~itself, all assertions follow, except for the final isomorphism claimed in part~ii).

For this, note that, for
every~$i$,
the subspaces
$\smash{H_{i,+}}$
and~$\smash{H_{i,-}}$
are both isotropic, while
$\smash{H_i = H_{i,+} \oplus H_{i,-}}$
is non-degenerate. Thus, the cup product pairing identifies 
$\smash{H_{i,-}}$
with the dual
$\smash{H_{i,+}^\vee}$.
But then, for an arbitrary element
$\smash{g \in \GL(H_{i,+})}$,
the map
$\smash{(g,(g^\vee)^{-1}) \in \GL(H_{i,+}) \times \GL(H_{i,+}^\vee) \subset \GL(H_i)}$
is orthogonal, and there is no other choice for the second component that would lead to this property.
\eop
}
\end{theo}

\begin{coro}
\label{ST0}
Let\/~$X$
be a\/
$K3$~surface
over\/~$\bbQ$,
$\smash{H \subset H^2(X(\bbC), \bbQ)}$
the transcendental part of the cohomology,
and\/~$E$
its endomorphism~field.

\begin{iii}
\item
If\/~$E$
is totally real of
degree\/~$\delta$~then
$\smash{(\ST^2(X))^0 \cong [\SO_{\frac{22-r}\delta}(\bbR)]^\delta}$.
For\/~$\delta=1$,
this includes the generic
case~$E=\bbQ.$
\item
If\/~$E$
is a CM~field of
degree\/~$2\delta$~then
$\smash{(\ST^2(X))^0 \cong [\U_{\frac{22-r}{2\delta}}]^\delta}$.
\end{iii}\smallskip

\noindent
{\bf Proof.}
{\em
The maximal compact subgroup of
$\SO_n(\bbC)$
is~$\SO_n(\bbR)$
and that of
$\GL_n(\bbC)$
is~$\U_n$,
cf.\ \cite[Table~(1.144)]{Kn}.
}
\eop
\end{coro}

\subsubsection*{Upper estimates for the component group}

\begin{lem}
\label{est_naive}
Let\/~$X$
be a\/
$K3$~surface
over\/~$\bbQ$,
$\smash{H_\tr \subset H^2(X(\bbC), \bbQ)}$
the transcendental part of the cohomology,
and\/~$E$
its endomorphism~field.

\begin{iii}
\item
If\/~$E$
is totally real~then\/
$N_{\O(H_\tr \otimes_\bbQ \bbC)} ((G_{X,l}^{2,\Zar}(\bbC))^0) = [\O(H_1) \times\cdots\times \O(H_\delta)] \rtimes S_\delta$,
the group\/
$S_\delta$
permuting the\/
$\delta$
direct~factors.
\item
If\/~$E$
is a CM~field~then
$$N_{\O(H_\tr \otimes_\bbQ \bbC)} ((G_{X,l}^{2,\Zar}(\bbC))^0) = [\O(H_1)_{(H_{1,+},H_{1,-})} \!\times\cdots\times \O(H_\delta)_{(H_{\delta,+},H_{\delta,-})}] \rtimes (\bbZ/2\bbZ)^\delta \!\rtimes S_\delta.$$
Here,
$e_i \in (\bbZ/2\bbZ)^\delta$
interchanges\/
$H_{i,+}$
with~$H_{i,-}$,
while\/
$S_\delta$
permutes the\/
$\delta$
direct~factors.
\end{iii}\smallskip

\noindent
{\bf Proof.}
{\em
i)
``$\supseteq$''
is clear.\smallskip

\noindent
``$\subseteq$'':
The natural action of the group
$\SO(H_1) \times\cdots\times \SO(H_\delta)$
on~$H_\tr \!\otimes_\bbQ\! \bbC$
setwise stabilises the subvector spaces
$H_1, \ldots, H_\delta$
and no others of
dimension~$\smash{\frac{22-r}\delta}$.
Hence, a linear map normalising
$\SO(H_1) \times\cdots\times \SO(H_\delta)$
must permute
$H_1, \ldots, H_\delta$.
As~every orthogonal map sending
$H_1, \ldots, H_\delta$
to themselves lies in
$\O(H_1) \times\cdots\times \O(H_\delta)$,
the assertion is~proven.\smallskip

\noindent
ii)
Again,
``$\supseteq$''
is clear.\smallskip

\noindent
``$\subseteq$'':
Here, the group
$\O(H_1)_{(H_{1,+},H_{1,-})} \times\cdots\times \O(H_\delta)_{(H_{\delta,+},H_{\delta,-})}$
stabilises the subvector spaces
$H_{1,+}, H_{1,-}, H_{2,+}, \ldots, H_{\delta,-}$
and no others of
dimension~$\smash{\frac{22-r}{2\delta}}$.
Thus, a linear map normalising
$\O(H_1)_{(H_{1,+},H_{1,-})} \times\cdots\times \O(H_\delta)_{(H_{\delta,+},H_{\delta,-})}$
must permute the spaces
$H_{1,+}, H_{1,-}, H_{2,+}, \ldots, H_{\delta,-}$.
Furthermore,~for
$i = 1,\ldots,\delta$,
the space
$H_{i,+}$
is perpendicular to both,
$H_{j,+}$
and~$H_{j,-}$,
for~$j\neq i$,
but it is not perpendicular
to~$H_{i,-}$.
Thus, the sets
$\{H_{i,+}, H_{i,-}\}$
form a block~system. The proof is therefore~complete.
}
\eop
\end{lem}

By construction, one has
$G_{X,l,\tr}^{2,\Zar}(\bbC) \subseteq \O(H_\tr \!\otimes_\bbQ\! \bbC)$.
Moreover, in every Lie group, the neutral component is a normal subgroup. Thus, in the RM as well as the CM cases, for the component group, one finds an~inclusion
\begin{align}
\label{comp_incl}
i_{X,l}\colon C_\tr := \ST^2_\tr(X)/(\ST^2(X))^0 \cong{} G_{X,l,\tr}^{2,\Zar}&(\bbC)/(G_{X,l}^{2,\Zar}(\bbC))^0 = G_{X,l,\tr}^{2,\Zar}/(G_{X,l}^{2,\Zar})^0\\[-1.5mm]
\hookrightarrow N_{\O(H_\tr \otimes_\bbQ \bbC)} ((&G_{X,l}^{2,\Zar}(\bbC))^0)/(G_{X,l}^{2,\Zar}(\bbC))^0 \cong (\bbZ/2\bbZ)^\delta \rtimes S_\delta \,. \nonumber
\end{align}
The~idea to simply compare with the normaliser is, of course, very~rough. In~fact, for
$\delta=2$
already,
$i_{X,l}(C_\tr) \subset (\bbZ/2\bbZ)^\delta \rtimes S_\delta$
is always a proper subgroup, as the next Theorem~shows.

\begin{rem}
We do not discuss in this article the question whether the homomorphism
$i_{X,l}$
is independent
of~$l$.
The~component
group~$C_\tr$
itself certainly is, as follows from \cite[p.~16, Th\'eor\`eme]{Se81}, cf.\ \cite[\S8.3.4]{Se12}.
\end{rem}

\begin{theo}
\label{upper_bound}
Let\/~$X$
be a\/
$K3$~surface
over\/~$\bbQ$.
Write\/~$E$
for the endomorphism field
of~$H_\tr \subset H^2(X(\bbC), \bbQ)$
and let\/
$E_0 \subseteq E$
be the maximal totally real subfield. Suppose~that
$E_0/\bbQ$
is Galois and that its Galois group is~cyclic. Then,~for any
prime\/~$l$
that is totally inert
in\/~$E_0$,
the following statements~hold.

\begin{abc}
\item
The~image of\/
$\pi_{X,l}\colon C_\tr \to S_\delta$
is a permutation group that is regular on any of its~orbits. In~other words, only the identity element has a fixed~point.
\item
Moreover,~the kernel of\/
$\pi_{X,l}$
is either trivial or of
order\/~$2$,
generated by the central element\/
$(-1,\dots,-1) \in (\bbZ/2\bbZ)^\delta$.
\end{abc}\smallskip

\noindent
{\bf Proof.}
{\em
a)
Suppose,~to the contrary, that there is an element
$\smash{A \in G_{X,l,\tr}^{2,\Zar}}$
that fixes an eigenspace
$H_i$,
but does not fix another,
$H_j$.
We~know that
$H_i$
and~$H_j$
are conjugate under
$\Gal(E/\bbQ)$.
As~$l$
is totally inert, this means
$\smash{\Frob_l^k(H_i) = H_j}$,
for a certain
$k\in\bbN$.
Since~$A$
is
\mbox{$\bbQ_l$-linear},
this~yields
$$H_j = \Frob_l^k(H_i) = \Frob_l^k(A(H_i)) = A(\Frob_l^k(H_i)) = A(H_j) \,,$$
a~contradiction.\smallskip

\noindent
b)
The kernel
of~$\pi$
consists of the elements stabilising each of the
$H_i$,
for
$i=1,\ldots,\delta$.
In~the CM case, the counter assumption is that some element
in~$\smash{G_{X,l,\tr}^{2,\Zar}}$
fixes the eigenspaces
$H_{i,+}$
and~$H_{i,-}$,
for
some~$i$,
but interchanges
$H_{j,+}$
and~$H_{j,-}$,
for a
certain~$j \neq i$.
This~is contradictory for exactly the same reason as in the proof of~a).

In~the RM case, the counter assumption is that there is some element
$\smash{A \in G_{X,l,\tr}^{2,\Zar}}$
being contained in
$\O(H_1) \times\cdots\times \O(H_\delta)$
and having
determinant~$1$
on some
$H_i$,
but
determinant~$(-1)$
on another,
$H_j$.
Again,~this is contradictory, as there is some
$k\in\bbN$
of the kind that
$\smash{\Frob_l^k(H_i) = H_j}$.

Indeed,~one has the
\mbox{$\bbQ_l$-linear}
map
$$\Exterior^\frac{22-r}\delta A\colon \Exterior^\frac{22-r}\delta(H_\tr \!\otimes_\bbQ\! \bbQ_l) \to \Exterior^\frac{22-r}\delta(H_\tr \!\otimes_\bbQ\! \bbQ_l) \,,$$
induced by
$A\colon H_\tr \!\otimes_\bbQ\! \bbQ_l \to H_\tr \!\otimes_\bbQ\! \bbQ_l$.
The base extension
to~$\overline\bbQ_l$
contains the one-dimensional subspaces
$\smash{\Exterior^\frac{22-r}\delta H_i}$,
on which 
$\smash{\Exterior^\frac{22-r}\delta A}$
acts as the identity, and 
$\smash{\Exterior^\frac{22-r}\delta H_j}$,
on which it acts as the multiplication
by~$(-1)$.
The~eigenspaces for the eigenvalues
$(+1)$
and
$(-1)$
are, however,
\mbox{$\bbQ_l$-subvector}
spaces
of~$\smash{\Exterior^\frac{22-r}\delta(H_\tr \!\otimes_\bbQ\! \bbQ_l)}$,
so that
$\smash{\Frob_l^k(\Exterior^\frac{22-r}\delta H_i) = \Exterior^\frac{22-r}\delta H_j}$
is impossible.
}
\eop
\end{theo}

\begin{exs}
\label{forbidden_elements}
\begin{iii}
\item
For~$\delta=2$,
one has that
$(\bbZ/2\bbZ)^2 \rtimes S_2$
is the dihedral group of order~eight. Part~b) of the Theorem forbids exactly two of its~elements. In~a somewhat symbolic notation, these are 
$\smash{\genfrac(){0pt}{}{- \,\,0}{\,0 \,\,+}}$
and
$\smash{\genfrac(){0pt}{}{+ \,\,0}{\,0 \,\,-}}$.
Thus, exactly two of the three conjugacy classes of subgroups of order four are still allowed, the cyclic subgroup being one of~them.

For~the actual occurrence of the cyclic group of order four, we refer to Example~\ref{ex5}, below. For~then other conjugacy class, there is a conjectural example presented as Example~\ref{ex6}.
\item
For~$\delta$
an odd prime, Theorem~\ref{upper_bound} implies that the component group is always cyclic of an order
dividing~$2\delta$.
\end{iii}
\end{exs}

\section{Experimental results}
\label{sec_exp}

\subsubsection*{The approach in general}

According to the Sato--Tate conjecture, one can use the theory of Lie groups in order to make a prediction on the distribution of the Frobenius traces. We~tested this in the situation of
$K3$~surfaces.
Depending on the Picard rank, the endomorphism field, and the jump character, various Lie groups occur, and hence various distributions are predicted. We~calculated the predicted densities as indicated in Section~\ref{Lie}.

To estimate the actual distributions, we used a Harvey style
\mbox{$p$-adic}
point counting algorithm~\cite{Ha} in order to determine the number of
$\bbF_{\!p}$-rational
points on the reduction
$\bmod\, p$,
for all primes
$p$
up to
$10^8$.
We~implemented the moving simplex idea~\cite[\S4.1]{Ha}, cf.\ \cite[Remark~4.8]{EJ16}.
In~order to speed up the computations, a
\mbox{$2$-adic}
algorithm was applied in addition~\cite{EJ22}. We~split the range
$[-6, 6]$
for the trace into 300 subintervals of equal lengths and counted the number of hits for each~subinterval. Representing the numbers of hits as columns, we then plotted the corresponding~histogram.

\subsubsection*{Running times}

It took around eight hours per surface on one core of an Intel(R) Core(TM)i7-7700 CPU processor running at
$3.6$\,GHz
to calculate the point counts in the case of a surface of type
$w^2 = xyz f_3(x,y,z)$.
For Example~\ref{ex5}, which is of the slightly more general shape
$w^2 = xy f_4(x,y,z)$,
it took 58~hours.

Note that the main step in the algorithm is to compute a small number of coefficients in huge powers
of~$f_6(x,y,z)$.
When working with a form of a particular shape as above, only the powers of a cubic, respectively quartic, form have to be considered, which leads to a massive reduction of the resulting computation.

\begin{rem}
For two of the seven
$K3$~surfaces
in our sample, the endomorphism fields are only conjectural. This~is not a serious problem, as this work is of a purely experimental character anyway. One~might consider the experiment as a test whether the correct Lie group is considered or whether blatant contradictions arise to the considerations above.
\end{rem}

\subsubsection*{Constraints concerning the endomorphism field}\leavevmode\medskip

\noindent
For~general considerations concerning the concept of the endomorphism field in the situation of a
$K3$~surface,
we refer to Paragraph~\ref{neut_comp}.

\begin{lem}
\label{e_fields}
Let\/~$X$
be a\/
$K3$~surface
over\/~$\bbC$.

\begin{abc}
\item
Suppose~that\/
$\rk\Pic X = 17$.
Then~the endomorphism field
is\/~$E=\bbQ$.
\item
Suppose~that\/
$\rk\Pic X = 16$.
Then~the endomorphism field
is either\/
$\bbQ$,
or a quadratic number field, or a CM~field of degree~six.
\end{abc}\smallskip

\noindent
{\bf Proof.}
{\em
One has
$\dim H_\tr = 22 - \rk\Pic X$.
Furthermore,
$[E:\bbQ] \mid \dim H_\tr$,
as
$H_\tr$
carries the structure of an
\mbox{$E$-vector}
space.\smallskip

\noindent
a)
Then
$[E:\bbQ] = 1$
or~$5$.
If~$[E:\bbQ] = 5$
then
$E$~is
not a CM~field, since
$5$
is~odd. Moreover,~in the RM~case, one has
$\smash{\frac{\dim H_\tr}{[E:\bbQ]} \geq 3}$~\cite[Lemma~3.2]{vG}.
Hence,~$E=\bbQ$.\smallskip

\noindent
b)
Then~$[E:\bbQ] = 1$,
$2$,
$3$,
or~$6$.
The~assumption
$[E:\bbQ] = 3$
is contradictory in exactly the same way as the assumption
$[E:\bbQ] = 5$
in~a). Moreover,~if
$[E:\bbQ] = 6$
then
$E$
cannot be totally real, since
$\smash{\frac{\dim H_\tr}{[E:\bbQ]} = 1 < 3}$.
}
\eop
\end{lem}

\begin{lem}
\label{CM_field}
Let\/~$X$
be a\/
$K3$~surface
over\/~$\bbC$.
Suppose~that\/
$X$
has CM by a quadratic
field\/~$\bbQ(\sqrt{-\delta})$,
for
$\delta \in \bbN$.

\begin{abc}
\item
If\/~$\dim H_\tr \equiv 2 \pmod 4$
then, for the discriminant\/~\cite[Chapitre~IV, \S1.1]{Se70}, one has\/
$\disc(H_\tr) = \overline\delta \in \bbQ^*/\bbQ^{*2}$.
\item
If\/~$\dim H_\tr \equiv 0 \pmod 4$
then\/
$\disc(H_\tr) = \overline1 \in \bbQ^*/\bbQ^{*2}$.
\end{abc}\smallskip

\noindent
{\bf Proof.}
{\em
Take an anisotopic vector
$v \in H_\tr$
and let
$I\colon H_\tr \to H_\tr$
be the endomorphism corresponding to
$\sqrt{-\delta}$.
Then
$(Iv,v) = (v,-Iv) = -(Iv,v)$,
i.e.\
$(Iv,v) = 0$,
by \cite[Theorem~1.5.1]{Za}. And similarly
$(Iv, Iv) = (v, -I^2v) = (v, \delta v) = \delta (v, v)$.
Thus,~the two-dimensional
\mbox{$I$-invariant}
quadratic subspace
$\langle v, Iv \rangle$
is of discriminant
$\smash{\big(\!\det \genfrac(){0pt}{}{(v, v)\;\;\; 0\;\;\;\;}{\;\;\;\;0 \;\;\delta(v, v)} \bmod \bbQ^{*2} \big) = \overline\delta \in \bbQ^*/\bbQ^{*2}}$.
The assertion follows inductively from~this.
}
\eop
\end{lem}

\subsubsection*{The surfaces inspected}
Each of the seven surfaces inspected is represented by a singular
degree~$2$
model of the shape
$$X'_i \colon w^2 = f_i(x,y,z) \,,$$
where
$f_i$,
for
$i=1,\ldots,7$,
is a ternary sextic form
over~$\bbQ$.
In~all cases, the ramification curve
$\smash{V(f_i) \subset \Pb^2_{\overline\bbQ}}$
has only ordinary double points. Thus, blowing up each of them once yields a
$K3$~surface~$X_i$
\cite[Theorem~8.2.27]{Do}, to which Lemma~\ref{bl_down}~applies.

Moreover, if there are
$N$
singular points then the exceptional curves
$E_1, \ldots, E_N$
together with the pull-back of a general line
in~$\smash{\Pb^2_\bbQ}$
generate a subgroup
of
rank~$(N+1)$
in~$\smash{\Pic X_{i,\overline\bbQ}}$.
If, in particular,
$V(f_i)$
geometrically splits into a union of six lines then
$\smash{\rk \Pic X_{i,\overline\bbQ} \geq 16}$.
If~$\smash{\rk \Pic X_{i,\overline\bbQ} = 16}$
holds exactly then
$$\disc \Pic X_{i,\overline\bbQ} = (\det\diag(2,-2,\ldots,-2) \bmod \bbQ^{*2}) = -\overline{1} \in \bbQ^*/\bbQ^{*2}$$
and hence
$\smash{\disc H_\tr = \overline{1} \in \bbQ^*/\bbQ^{*2}}$.
Thus,
$\smash{\bbQ(\sqrt{-1})}$
is the only imaginary quadratic field that is possible for~CM.

We~list the bad primes as well as the jump character for each of the seven sample surfaces in a table at the very end of this article. By bad primes, those of the obvious model
over~$\bbZ$
are meant, which is constructed from the double cover
of~$\Pb^2_\bbZ$,
defined by the equation
$f_i=0$,
by the blow-ups centred in the Zariski closures of the finitely many singular points of the generic~fibre. The~jump characters are obtained using \cite[Algorithm~2.6.1]{CEJ}. Note~that, in each case, not only the geometric Picard rank is known, but the geometric Picard group as a
\mbox{$\smash{\Gal(\overline\bbQ/\bbQ)}$-module}.\smallskip\pagebreak[3]

\subsubsection*{A generic example of Picard rank 16}

\begin{ex}
\label{ex1}
Let~$X_1'$
be the double cover
of~$\Pb^2_\bbQ$,
given by
$$w^2 = xyz (x+y+z) (3x+5y+7z) (-5x+11y-2z)$$
and
$X_1$
the
$K3$~surface
obtained as the minimal desingularisation
of~$X_1'$.

\begin{abc}
\item
Then the geometric Picard rank
of~$X_1$
is~$16$.
\item
The endomorphism field
of~$X_1$
is~$E=\bbQ$.
\end{abc}\smallskip

\noindent
{\bf Proof.}
a)
One~has a lower bound
of~$16$,
as the ramification locus has 15 singular~points. An~upper bound
of~$16$
is provided by the reduction
modulo~$31$,
which is of geometric Picard
rank~$16$.\smallskip

\noindent
b)
The reduction
modulo~$17$
is of geometric Picard
rank~$18$,
which, by \cite[Lemma~6.2]{EJ20a} implies that
$[E:\bbQ] \leq 2$.
Furthermore,~RM is excluded, since there is a reduction of geometric Picard rank
$16$~\cite[Corollary~4.12]{EJ14}. Finally,~if
$E$~were
a CM field then, by Lemma~\ref{CM_field}, the only option would be
$\smash{E = \bbQ(\sqrt{-1})}$.

In that case, one would have a decomposition
$H_\tr \!\otimes_\bbQ\! \bbC = H_+ \oplus H_-$,
the summands being eigenspaces for the eigenvalues
$\smash{\pm\sqrt{-1}}$,
and hence defined
over~$\smash{\bbQ(\sqrt{-1})}$.
By~Lemma~\ref{est_naive}, the algebraic monodromy group
$\smash{G_{X,l,\tr}^{2,\Zar}}$
has at most two~components. The non-neutral component, if present, interchanges the eigenspaces and hence all elements are of trace zero. The neutral component stabilises
$H_+$
and~$H_-$
and hence, the characteristic polynomial of every element factors
over~$\smash{\bbQ_l(\sqrt{-1})}$
into two cubic~polynomials.
However, the characteristic polynomial of
$\smash{\varrho^2_{X,17,\tr}(\Frob_{31}) \in G_{X,17,\tr}^{2,\Zar}}$
is
$\smash{t^6 - \frac{10}{31}t^5 + \frac1{31}t^4 + \frac{20}{31}t^3 + \frac1{31}t^2 - \frac{10}{31}t + 1}$,
which splits over
$\bbQ_{17}(\sqrt{-1}) = \bbQ_{17}$
into irreducible factors of degrees two and~four. Moreover, the trace of
$\smash{\varrho^2_{X,17,\tr}(\Frob_{31})}$
is nonzero, a contradiction.
\eop
\end{ex}

In view of the results above, Corollary~\ref{ST0} shows that
$(\ST^2(X_1))^0 \cong \SO_6(\bbR)$.
Moreover,
$\ST^2_\tr(X_1)/(\ST^2(X_1))^0 = \bbZ/2\bbZ$,
i.e.\
$\ST^2_\tr(X_1) \cong \O_6(\bbR)$.
Indeed,~for the component group, we have an upper bound
of~$\bbZ/2\bbZ$
by~(\ref{comp_incl}), and the trivial group is excluded, due to the nontrivial jump character, cf.\ Table~\ref{badpr_jump}.

\subsubsection*{An example of Picard rank 16 with trivial jump character}

\begin{ex}
\label{ex2}
Let~$X_2'$
be the double cover
of~$\Pb^2_\bbQ$,
given by
$$w^2 = xyz (2x+4y-3z) (x-5y-3z) (x+3y+3z)$$
and
$X_2$
the
$K3$~surface
obtained as the minimal desingularisation
of~$X_2'$.

\begin{abc}
\item
Then the geometric Picard rank
of~$X_2$
is~$16$.
\item
The endomorphism field
of~$X_2$
is~$E=\bbQ$.
\end{abc}\smallskip

\noindent
{\bf Proof.}
a)
One~has a lower bound
of~$16$,
as the ramification locus has 15 singular~points. An~upper bound
of~$16$
is provided by the reduction
modulo~$19$,
which is of geometric Picard
rank~$16$.\smallskip

\noindent
b)
The reduction
modulo~$13$
is of geometric Picard
rank~$18$,
which, as in Example~\ref{ex1}, leaves
$\smash{E = \bbQ(\sqrt{-1})}$
as the only nontrivial~option. Moreover,~this is excluded by observing that the characteristic polynomial of
$\smash{\varrho^2_{X_2,13,\tr}(\Frob_{19}) \in G_{X_2,13,\tr}^{2,\Zar}}$
is
$\smash{t^6 + \frac2{19}t^5 + \frac{13}{19}t^4 - \frac4{19}t^3 + \frac{13}{19}t^2 + \frac2{19}t + 1}$,
which splits over
$\bbQ_{13}(\sqrt{-1}) = \bbQ_{13}$
into irreducible factors of degrees two and four. Note~that the trace is nonzero.
\eop
\end{ex}

Here,~one has
$\smash{\ST^2_\tr(X_2) \cong \SO_6(\bbR)}$.
Indeed, as above,
$\smash{(\ST^2(X_2))^0 \cong \SO_6(\bbR)}$,
and the component group is trivial, due to the trivial jump character, cf.\ Table~\ref{badpr_jump}.

\subsubsection*{An example of Picard rank 17 with trivial jump character}

\begin{ex}
\label{ex3}
Let~$X_3'$
be the double cover
of~$\Pb^2_\bbQ$,
given by
$$w^2 = xyz (4x+9y+z)(-x-y-4z)(16x+25y+z)$$
and
$X_3$
the
$K3$~surface
obtained as the minimal desingularisation
of~$X_3'$.

\begin{abc}
\item
Then the geometric Picard rank
of~$X_3$
is~$17$.
\item
The endomorphism field
of~$X_3$
is~$E=\bbQ$.
\end{abc}\smallskip

\noindent
{\bf Proof.}
a)
The 16 elements
$\smash{\pi^*[l], [E_1], \ldots, [E_{15}] \in \Pic X_{3,\overline\bbQ}}$
are linearly independent, as~before. Moreover, the inverse image of the conic
$C := \Vb(xy + yz + zw) \subset \Pb^2$~in~$X_3$
splits into two curves,
$\smash{C'}$
and~$\smash{C''}$,
as a Gr\"obner base calculation~shows.

We claim that
$[C']$
is independent of the 16 elements~above. Indeed, otherwise
$[C']$
would be invariant under the involution of the double
cover~$\pi$.
Since~one has
$[C'] + [C''] = 2\pi^*[l]$
and
$C'$
is interchanged with
$C''$
under the involution, this implies
$\smash{[C'] = \pi^*[l] \in \Pic (X_3)_{\overline\bbQ}}$.
But
$C'$
is rational, and hence 
a~$(-2)$-curve,
while
$\pi^*[l]$
has self-intersection
number~$(+2)$,
a~contradiction. Thus, there is a lower bound
of~$17$.

Concerning the upper bound, the reductions modulo
$13$
and~$23$
are both of geometric Picard
rank~$18$.
The~characteristic polynomials of the Frobenii are
$$\textstyle (t-1)^{18} (t^4 + \frac{20}{13}t^3 + \frac{30}{13}t^2 + \frac{20}{13}t + 1)
\quad{\rm and}\quad
(t-1)^{18} (t^4 + \frac{36}{23}t^3 + \frac{42}{23}t^2 + \frac{36}{23}t + 1) \,,$$
so that the Artin--Tate formula~\cite[Theorem~6.1]{Mi} determines the discriminants of the four-dimensional lattices to
$\overline{6}$
and
$\overline{10} \in \bbQ^*/\bbQ^{*2}$,
respectively. I.e., the lattices are incompatible and van Luijk's method~\cite{vL} lets the upper bound drop
to~$17$.\smallskip

\noindent
b) follows immediately from~a), in view of Lemma~\ref{e_fields}.a).
\eop
\end{ex}

Corollary~\ref{ST0} shows that
$(\ST^2(X_3))^0 \cong \SO_5(\bbR)$.
This~is the only component, as the jump character is trivial, cf.\ Table~\ref{badpr_jump}.\vspace{2.5mm}\smallskip\pagebreak[3]

\subsubsection*{The generic trace distributions}
~\vspace{-4.2mm}

\begin{figure}[H]
\begin{center}
\includegraphics[scale=0.47]{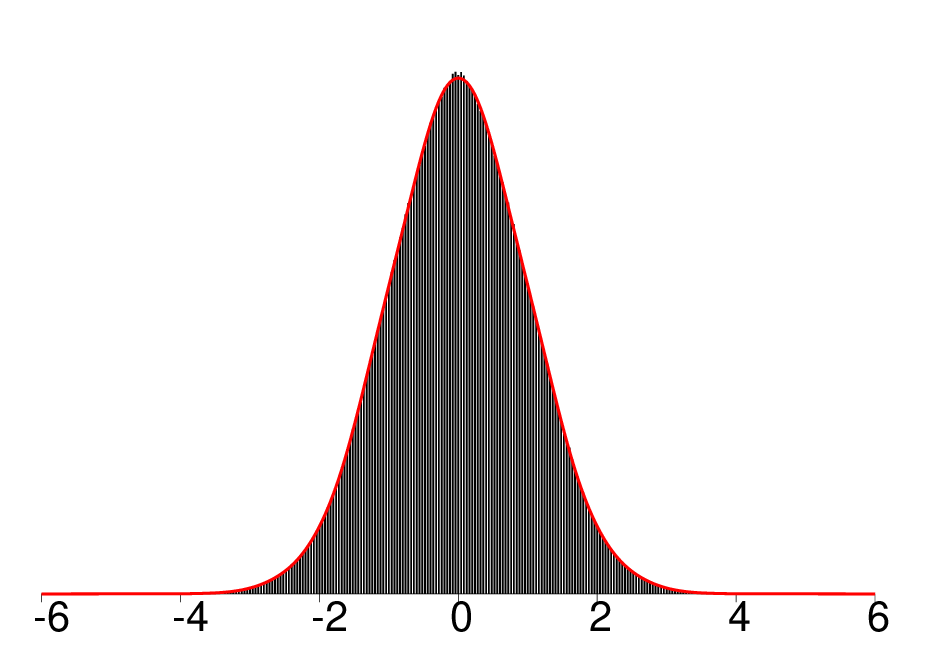}
\end{center}
\end{figure}\vspace{-9.2mm}

\begin{figure}[H]
\begin{center}
\includegraphics[scale=0.47]{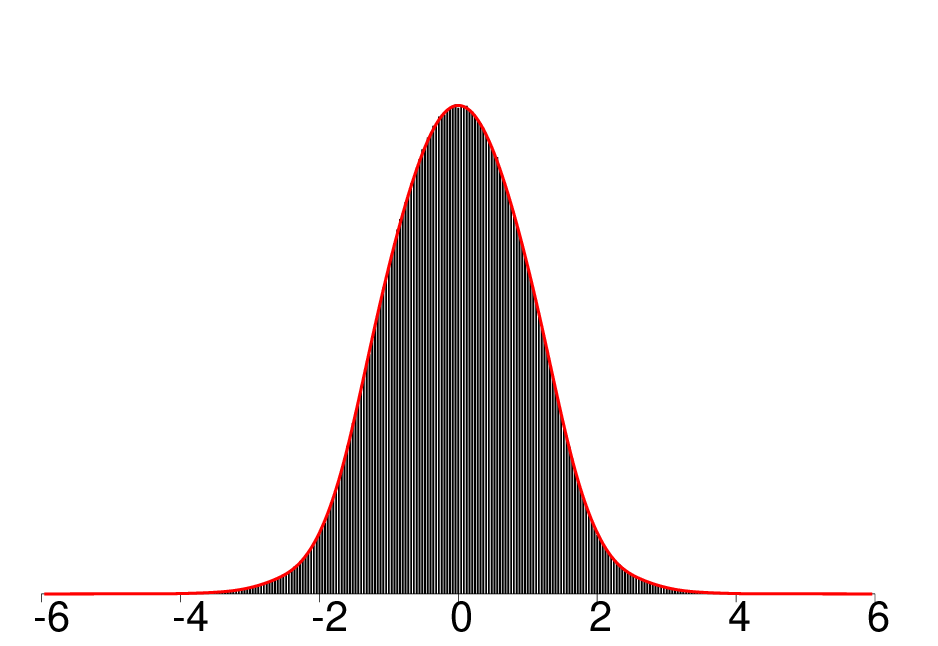}
\end{center}
\end{figure}\vspace{-9.2mm}

\begin{figure}[H]
\begin{center}
\includegraphics[scale=0.47]{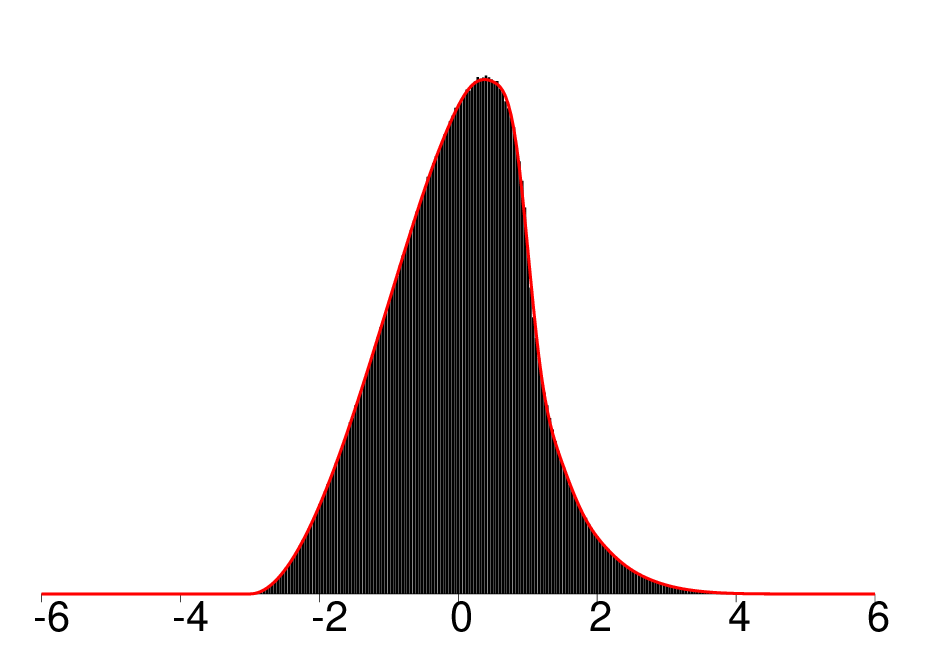}
\end{center}\vspace{-5mm}
\caption{Trace distributions for Examples~\ref{ex1}, \ref{ex2}, and \ref{ex3}}\label{generisch_plot}
\end{figure}\vspace{-2mm}

The red lines in Figure~\ref{generisch_plot} show the densities of the theoretical trace distributions, as predicted by the Sato--Tate conjecture. For example~\ref{ex1}, one has the superposition of the distributions for the two components, as explained in Section~\ref{Lie}. Note~that the theoretical density for Example~\ref{ex3} is not symmetric.

\subsubsection*{An example with CM by\/
$\bbQ(\sqrt{-1})$}

\begin{ex}
\label{ex4}
Let~$X_4'$
be the double cover
of~$\Pb^2_\bbQ$,
given by
$$w^2 = xyz (x+y+z) (x+2y+3z) (5x+8y+20z)$$
and
$X_4$
the
$K3$~surface
obtained as the minimal desingularisation
of~$X_4'$.

\begin{abc}
\item
Then the geometric Picard rank
of~$X_4$
is~$16$.
\item
The endomorphism field
of~$X_4$
is~$E=\bbQ(\sqrt{-1})$.
\end{abc}\smallskip

\noindent
{\bf Proof.}
a)
One~has a lower bound
of~$16$,
as the ramification locus has 15 singular~points. An~upper bound
of~$16$
is provided by the reduction
modulo~$13$,
which is of geometric Picard
rank~$16$.\smallskip

\noindent
b)
The reduction
modulo~$11$
is of geometric Picard
rank~$18$,
which, as in Example~\ref{ex1}, leaves
$\smash{E = \bbQ(\sqrt{-1})}$
as the only nontrivial~option. The~result follows from Theorem~\ref{Qi_deg2} below, as the endomorphism field does not shrink under specialisation~\cite[Corollary~4.6]{EJ20a}.
\eop
\end{ex}

Here, Corollary~\ref{ST0} yields
$(\ST^2(X_4))^0 \cong \U_3$.
Moreover, for the component group, one has
$\ST^2_\tr(X_4)/(\ST^2(X_1))^0 = \bbZ/2\bbZ$.
Indeed,~(\ref{comp_incl}) gives an upper bound
of~$\bbZ/2\bbZ$,
and there must be a second component, due to the nontrivial jump character, cf.\ Table~\ref{badpr_jump}.\vspace{-3.5mm}

\begin{figure}[H]
\begin{center}
\includegraphics[scale=0.47]{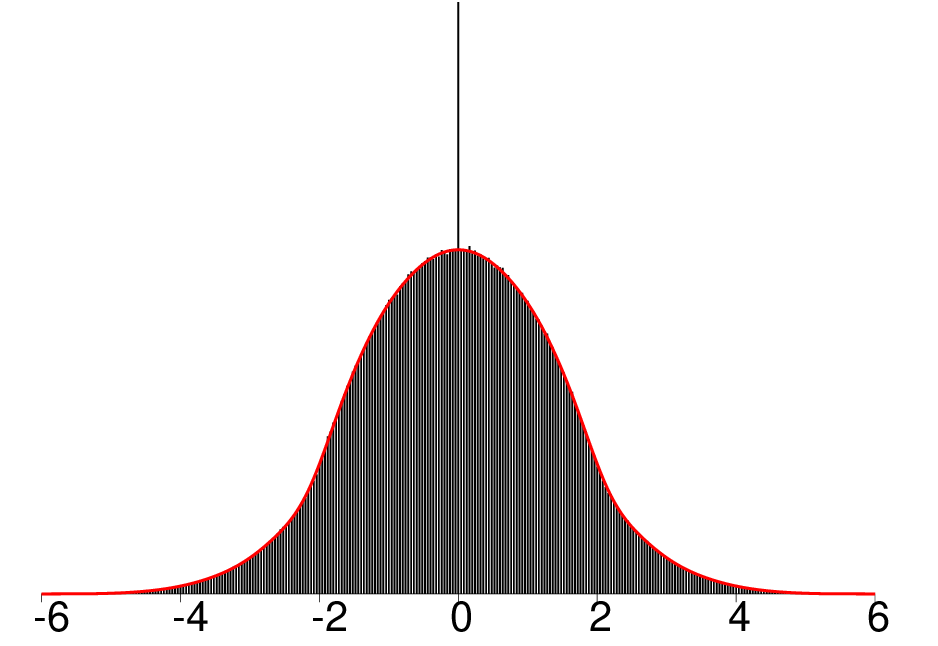}
\end{center}\vspace{-5mm}
\caption{Trace distribution for Example~\ref{ex4}}\label{CM_sqrt_m1_plot}
\end{figure}\vspace{-3mm}

In the figure above, the spike is of
mass~$1/2$.

\subsubsection*{An example with RM and a cyclic component group of order four}

\begin{ex}
\label{ex5}
Let~$X_5'$
be the double cover
of~$\Pb^2_\bbQ$,
given by
\begin{align*}
w^2 = xy (x^4 - 7x^3y - x^3z + 19x^2y^2 + 4x^2yz + x^2z^2 - 23xy^3 - 7xy^2z - 6xyz^2& \\[-1.5mm]
{}- xz^3 + 11y^4 + 7y^3z + 9y^2z^2 + 3yz^3 &+ z^4)
\end{align*}
and
$X_5$
the
$K3$~surface
obtained as the minimal desingularisation
of~$X_5'$.

\begin{abc}
\item
Then the geometric Picard rank
of~$X_5$
is~$16$.
\item
The endomorphism field
of~$X_5$
is~$E=\bbQ(\sqrt{5})$.
\end{abc}\smallskip

\noindent
{\bf Proof.}
a)
One~has a lower bound
of~$16$,
as the ramification locus has 15 singular~points.
As far as upper bounds are concerned, the reductions modulo
$7$
and~$19$
are both of geometric Picard
rank~$18$.
The~characteristic polynomials of the Frobenii are
$$\textstyle (t-1)^6 (t+1)^4 (t^2+1)^4 (t^4 - \frac{12}7 t^2 + 1)
\quad{\rm and}\quad
(t-1)^{10} (t+1)^8 (t^4 - \frac4{19} t^3 + \frac{34}{19} t^2 - \frac4{19} t + 1) \,,$$
so that the Artin--Tate formula~\cite[Theorem~6.1]{Mi} determines the discriminants of the four-dimensional lattices to
$\overline{1}$
and
$\overline{5} \in \bbQ^*/\bbQ^{*2}$,
respectively. I.e., the lattices are incompatible and van Luijk's method~\cite{vL} lets the upper bound drop
to~$17$.

On~the other hand, the endomorphism field
of~$X_5$
contains
$\smash{\bbQ(\sqrt{5})}$,
which excludes the option of
rank~$17$.
Indeed,
$X_5$
is isomorphic to the specialisation
to~$t=0$
of the family described in \cite[Example~1.5]{EJ20a}. The~isomorphism is induced by the automorphism
of~$\Pb^2_{\overline\bbQ}$,
given by the~matrix
$$
\left(
\begin{array}{rrr}
1 & \;2 & \;0 \\
0 &   1 &   0 \\
0 &   0 &   1
\end{array}
\right) \,.\medskip
$$

\noindent
b)
In~view of~a), this follows from \cite[Example~1.5.iv)]{EJ20a}.
\eop
\end{ex}

Corollary~\ref{ST0} shows that
$(\ST^2(X_5))^0 \cong [\SO_3(\bbR)]^2$.

\begin{theo}
The representation\/
$\smash{\varrho^2_{X_5,l,\tr}}$
induces an~isomorphism
$$\overline\varrho\colon \Gal(\bbQ(\zeta_5)/\bbQ) \longrightarrow \ST^2_\tr(X_5)/(\ST^2(X_5))^0 \,.$$

\noindent
{\bf Proof.}
{\em
{\em First step.}
Generalities.

\noindent
The jump character is trivial, so
$\ST^2_\tr(X_5)/(\ST^2(X_5))^0$
is bound to the cosets
$\smash{\genfrac(){0pt}{}{+ \,\,0}{\,0 \,\,+}}$,
$\smash{\genfrac(){0pt}{}{- \,\,0}{\,0 \,\,-}}$,
$\smash{\genfrac(){0pt}{}{\,0 \,\,+}{- \,\,0}}$,
and
$\smash{\genfrac(){0pt}{}{\,0 \,\,-}{+ \,\,0}}$.
Quite generally, there exists a unique number field
$L_0$,
for which
$\smash{\varrho^2_{X_5,l,\tr}}$
induces an isomorphism
$\smash{\overline\varrho\colon \Gal(L_0/\bbQ) \stackrel{\cong}{\longrightarrow} \ST^2_\tr(X_5)/(\ST^2(X_5))^0}$,
cf.\ Remark~\ref{rems_allg}.i). In~our situation, we find that
$\smash{L_0}$
is cyclic of a degree dividing~four.\smallskip

\noindent
{\em Second step.}
$\smash{L_0 \supsetneqq \bbQ(\sqrt{5})}$.

\noindent
According to Chebotarev, the elements
$\smash{\tau^{-1}\Frob_p\tau}$,
for
$\smash{\tau\in\Gal(\overline\bbQ/\bbQ)}$
and
$p \equiv 2,3 \pmod 5$,
are dense in the nontrivial coset
$\smash{C_5 := \Gal(\overline\bbQ/\bbQ) \setminus U_5}$
of the open subgroup
$\smash{U_5 := \Gal(\overline\bbQ/\bbQ(\sqrt{5})) \subset \Gal(\overline\bbQ/\bbQ)}$
of index~two. Moreover,~\cite[Lemma~6.7]{EJ20a} shows together with Lemma~\ref{bl_down} that
$\smash{\Tr(\varrho^2_{X_5,l,\tr}(\Frob_p)) = 0}$,
for every prime
$p \equiv 2,3 \pmod 5$.
Consequently,
\begin{equation}
\label{Tr0}
\Tr(\varrho^2_{X_5,l,\tr}(\sigma)) = 0 \,,
\end{equation}
for every
$\sigma \in C_5$.

Since,
$\smash{U_5 \subset \Gal(\overline\bbQ/\bbQ)}$
is a subgroup of finite index,
$\smash{\overline{\varrho^2_{X_5,l,\tr}(U_5)} \subseteq G_{X_5,l,\tr}^{2,\Zar}}$
has the same neutral component, only the component group may differ. Moreover, due to~(\ref{Tr0}),
$\smash{\overline{\varrho^2_{X_5,l,\tr}(C_5)} \subseteq G_{X_5,l,\tr}^{2,\Zar}}$
is certainly a nontrivial coset. In particular,
$\smash{\overline{\varrho^2_{X_5,l,\tr}(U_5)}}$
must be a proper subgroup of
$\smash{G_{X_5,l,\tr}^{2,\Zar}}$,
which yields that
$\smash{L_0 \supseteq \bbQ(\sqrt{5})}$.

Furthermore, (\ref{Tr0}) shows that the coset
$\smash{\overline{\varrho^2_{X_5,l,\tr}(C_5)}}$
consists only of components of type
$\smash{\genfrac(){0pt}{}{\,0 \,\,+}{- \,\,0}}$
and
$\smash{\genfrac(){0pt}{}{\,0 \,\,-}{+ \,\,0}}$.
In particular,
$\ST^2_\tr(X_5)/(\ST^2(X_5))^0$
is indeed of order~four.\smallskip

\noindent
{\em Third step.}
Conclusion.

\noindent
A~standard argument involving the smooth specialisation theorem for \'etale cohomology groups \cite[Exp.~XVI, Corollaire~2.2]{SGA4} shows that
$L_0$
is unramified at every prime
$p \neq 2$,
$5$,
and~$l$,
cf.\ \cite[Lemma~2.2.3.a)]{CEJ}. The~field
$L_0$
is, moreover, known to be independent
of~$l$~\mbox{\cite[p.~16, Th\'eor\`eme]{Se81}},
cf.\ \cite[\S8.3.4]{Se12}. Thus, working with
$l=2$
or~$5$,
one finds that
$L_0$
may ramify only at
$2$
and~$5$.

Besides
$\smash{\bbQ(\zeta_5)}$,
there are only three cyclic number fields of degree four that are unramified outside
$2$
and~$5$
and
contain~$\smash{\bbQ(\sqrt{5})}$.
These are the quadratic twists of
$\smash{\bbQ(\zeta_5)}$
by
$\smash{\bbQ(\sqrt\delta)}$,
for
$\delta = -1$,
$2$,
and~$(-2)$.
I.e., the unique further cyclic subfield of degree four
in~$\smash{\bbQ(\zeta_5, \sqrt\delta)}$.
Indeed, let
$L$
be such a field. Then, since
$\smash{\sqrt{5} \in L}$
and
$\smash{\sqrt{5} \in \bbQ(\zeta_5)}$,
the field
$\smash{L(\zeta_5)}$
has Galois group
$\bbZ/4\bbZ \times \bbZ/2\bbZ$.
Thus,
$\smash{L(\zeta_5) = \bbQ(\zeta_5, \sqrt\delta)}$,
for some
$\delta\in\bbZ$.
The claim follows, as
$\smash{L(\zeta_5)}$
is unramified outside
$2$
and~$5$.

Suppose that
$L_0$
is the quadratic twist of
$\smash{\bbQ(\zeta_5)}$
by
$\smash{\bbQ(\sqrt\delta)}$,
for
$\delta = -1$,
$2$,
or~$(-2)$.
Then
$\Frob_p \in \Gal(L_0/\bbQ)$
is not the neutral element for
$p=11$
in the first two cases, and for
$p=31$
in the third. However, an experiment shows that
$\smash{\varrho^2_{X_5,l,\tr}(\Frob_{11})}$
and
$\smash{\varrho^2_{X_5,l,\tr}(\Frob_{31})}$
are contained in the neutral component, which completes the~proof.
}
\eop
\end{theo}

\subsubsection*{An example with RM and the Klein four group as the component group}

\begin{ex}
\label{ex6}
Let~$X_6'$
be the double cover
of~$\Pb^2_\bbQ$,
given by
$$w^2 = xyz (x^3 - 14x^2z + 11xy^2 - xz^2 + 12y^3 - 14y^2z - 12yz^2 + 14z^3)$$
and
$X_6$
the
$K3$~surface
obtained as the minimal desingularisation
of~$X_6'$.

\begin{abc}
\item
Then the geometric Picard rank
of~$X_6$
is~$16$.
\item
The endomorphism field
of~$X_6$
is at most~quadratic.
\end{abc}\smallskip

\noindent
{\bf Proof.}
a)
For~the lower bound, the situation is analogous to \cite[Example~2.7.3]{CEJ}. One~immediately has a lower bound
of~$13$,
as the ramification locus has twelve singular~points. Among~them, ten are
\mbox{$\bbQ$-rational},
the two others are defined over
$\smash{\bbQ(\sqrt{-47})}$,
and conjugate to each~other. Moreover,~there are a
\mbox{$\bbQ$-rational}
line, the inverse image of which splits over
$\smash{\bbQ(\sqrt{14})}$,
and two conics that are defined over
$\smash{\bbQ(\sqrt{-1})}$
and conjugate to each other, the inverse images of which split over
$\smash{\bbQ(\sqrt{-1}, \sqrt{42})}$.
Thus, there is a sublattice
$\smash{P \subseteq \Pic X_{6,\overline\bbQ}}$
of
rank~$16$,
such~that
$\smash{P \otimes_\bbZ \bbC = \chi_\triv^{12} \oplus \chi_{\bbQ(\sqrt{-47})} \oplus \chi_{\bbQ(\sqrt{14})} \oplus \chi_{\bbQ(\sqrt{42})} \oplus \chi_{\bbQ(\sqrt{-42})}}$
(cf.\ \cite{CEJ} for~notation).
It~is a routine work that was carried out with some help of the machine to set up an intersection matrix and to calculate
that~$\disc(P \!\otimes_\bbZ\! \bbQ) = (-3) \in \bbQ^*/\bbQ^*{}^2$.

Concerning~an upper bound, the reductions modulo
$19$
and~$59$
are both of geometric Picard
rank~$18$.
The~characteristic polynomials of the Frobenii are
$$\textstyle (t-1)^{14} (t+1)^4 (t^4 + \frac{36}{19}t^2 + 1)
\quad{\rm and}\quad
(t-1)^{14} (t+1)^4 (t^4 - \frac{116}{59}t^3 + \frac{162}{59}t^2 - \frac{116}{59}t + 1) \,,$$
so that the Artin--Tate formula~\cite[Theorem~6.1]{Mi} determines the discriminants of the four-dimensional lattices to
$\overline{1}$
and
$\overline{6} \in \bbQ^*/\bbQ^{*2}$,
respectively. I.e., the lattices are incompatible and van Luijk's method~\cite{vL} lets the upper bound drop
to~$17$.

At~this point, a modification of the method described in \cite{EJ11} allows to reduce the upper bound even~further. For~this, suppose that one had
$\smash{\rk \Pic X_{6,\overline\bbQ} = 17}$.
The~$\smash{\Gal(\overline\bbQ/\bbQ)}$-representation
$\smash{\Pic X_{6,\overline\bbQ} \otimes_\bbZ \bbQ}$
then splits off a one-dimensional direct summand
$\smash{V \subset (P \otimes_\bbZ \bbQ)^\perp}$.

Let~us particularly consider the action
of~$\smash{\Frob_{19} \in \Gal(\overline\bbQ/\bbQ)}$.
The~characteristic polynomial on the whole
of~$\smash{\Pic X_{6,\overline\bbQ} \otimes_\bbZ \bbQ}$
is
then~$\smash{(t-1)^{14} (t+1)^4 (t^4 + \frac{36}{19}t^2 + 1)}$.
Furthermore, as
$(-47)$,
$14$,
$(-1)$,
and~$(-42)$
are all quadratic non-residues
modulo~$19$,
the action splits
$\smash{P \otimes_\bbZ \bbQ}$
into a
\mbox{$13$-dimensional}
invariant subspace
$(P \otimes_\bbZ \bbQ)^+$
and a three-dimensional
$(-1)$-eigenspace
$(P \!\otimes_\bbZ\! \bbQ)^-$.
Having set up in {\tt magma} the corresponding intersection matrices with respect to suitable bases, one calculates that
$\smash{\disc((P \!\otimes_\bbZ\! \bbQ)^+) = \overline{2} \in \bbQ^*/\bbQ^*{}^2}$
and
$\smash{\disc((P \!\otimes_\bbZ\! \bbQ)^-) = (-\overline{6}) \in \bbQ^*/\bbQ^*{}^2}$.

Moreover,~the characteristic polynomial
of~$\Frob_{19}$
on
$\smash{(P \otimes_\bbZ \bbQ)^\perp}$
turns out to be
$\smash{(t-1) (t+1) (t^4 + \frac{36}{19}t^2 + 1)}$.
Therefore,~$V$
may only be one of the one-di\-men\-sional eigenspaces, either
$V^+_{19}$
or~$V^-_{19}$.
On~the other hand, an application of the Artin--Tate formula shows that
$\smash{\disc((P \!\otimes_\bbZ\! \bbQ)^+ \!\!\perp\! V^+_{19}) = (-\overline{74}) \in \bbQ^*/\bbQ^*{}^2}$.
I.e.,~that
$\smash{\disc(V^+_{19}) = (-\overline{37}) \in \bbQ^*/\bbQ^*{}^2}$.
Finally, the Artin--Tate formula
for~$\bbF_{\!19^2}$
yields
$\smash{\disc((P \!\otimes_\bbZ\! \bbQ) \!\perp\! V^+_{19} \!\perp\! V^-_{19}) = (-\overline{1}) \in \bbQ^*/\bbQ^*{}^2}$,
so that
$\smash{\disc(V^-_{19}) = (-\overline{111}) \in \bbQ^*/\bbQ^*{}^2}$
results. Consequently,
$\smash{\disc(V) = (-\overline{37})}$
or~$\smash{(-\overline{111}) \in \bbQ^*/\bbQ^*{}^2}$.

From~this, a contradiction arises when one repeats the argument for a suitable second prime~number. For~example, the action
of~$\Frob_{127}$
on~$\smash{P \!\otimes_\bbZ\! \bbQ}$
has exactly the same invariant subspace
$(P \!\otimes_\bbZ\! \bbQ)^+$.
Moreover,~on
$\smash{(P \!\otimes_\bbZ\! \bbQ)^\perp}$,
both the
$(+1)$-
and
\mbox{$(-1)$-eigenspaces}
are again of dimension~one. A~calculation completely analogous to the one above indicates that nothing but
$\smash{\disc(V) = (-\overline{229})}$
or~$(-\overline{687}) \in \bbQ^*/\bbQ^*{}^2$~may
happen. This~provides the desired contradiction and hence completes the proof of~a).\smallskip

\noindent
b)
As~there are reductions of
rank~$18$,
this is a consequence of \cite[Lemma 6.2]{EJ20a}.
\eop
\end{ex}

There is strong evidence that the endomorphism field of
$X_6$
is in fact
$E=\bbQ(\sqrt{3})$.
The evidence has been described in \cite[Section~5]{EJ16}. Note~that
$\smash{X_6 = V^{(3)}_{1,2}}$
in the notation of \cite[Conjectures~5.2]{EJ16}. Thus, conjecturally,
$(\ST^2(X_6))^0 \cong [\SO_3(\bbR)]^2$.

The~observation that
$\smash{\Tr(\varrho^2_{X_6,l,\tr}(\Frob_p)) = 0}$
for all primes
$p = \pm5 \pmod {12}$
has meanwhile been extended
to~$p < 10^8$.
As~these are exactly the primes, at which the jump character evaluates
to~$(-1)$,
the component group
$\ST^2_\tr(X_6)/(\ST^2(X_6))^0$
is bound to the elements written symbolically as
$\smash{\genfrac(){0pt}{}{+ \,\,0}{\,0 \,\,+}}$,
$\smash{\genfrac(){0pt}{}{- \,\,0}{\,0 \,\,-}}$,
$\smash{\genfrac(){0pt}{}{\,0 \,\,+}{+ \,\,0}}$,
and
$\smash{\genfrac(){0pt}{}{\,0 \,\,-}{- \,\,0}}$.

The~component
$\smash{\genfrac(){0pt}{}{- \,\,0}{\,0 \,\,-}}$
is indeed met, thus the component group is isomorphic to the Klein four~group. According to our experiments,
$\smash{\overline{x}_p \in [\O_3^-(\bbR)]^2}$
if and only if
$p \equiv \pm 1 \pmod {12}$
and
$\smash{(\frac{-2\cdot7\cdot47}p) = -1}$.\vspace{2.5mm}

\subsubsection*{The trace distributions in the RM examples}
~\vspace{-1.0mm}

\begin{figure}[H]
\begin{center}
\includegraphics[scale=0.47]{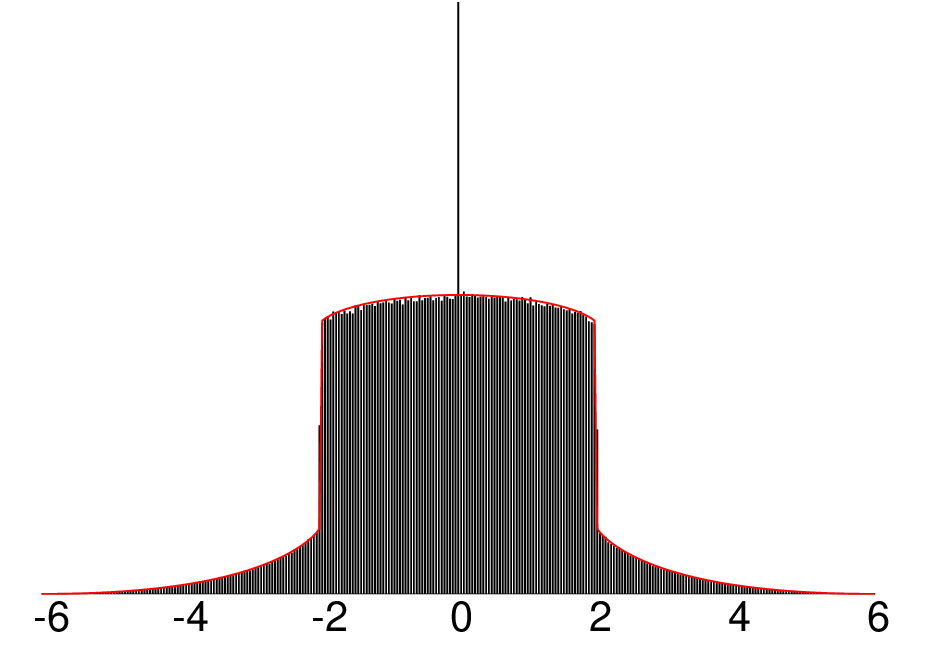}
\end{center}
\end{figure}\vspace{-4.7mm}

\begin{figure}[H]
\begin{center}
\includegraphics[scale=0.47]{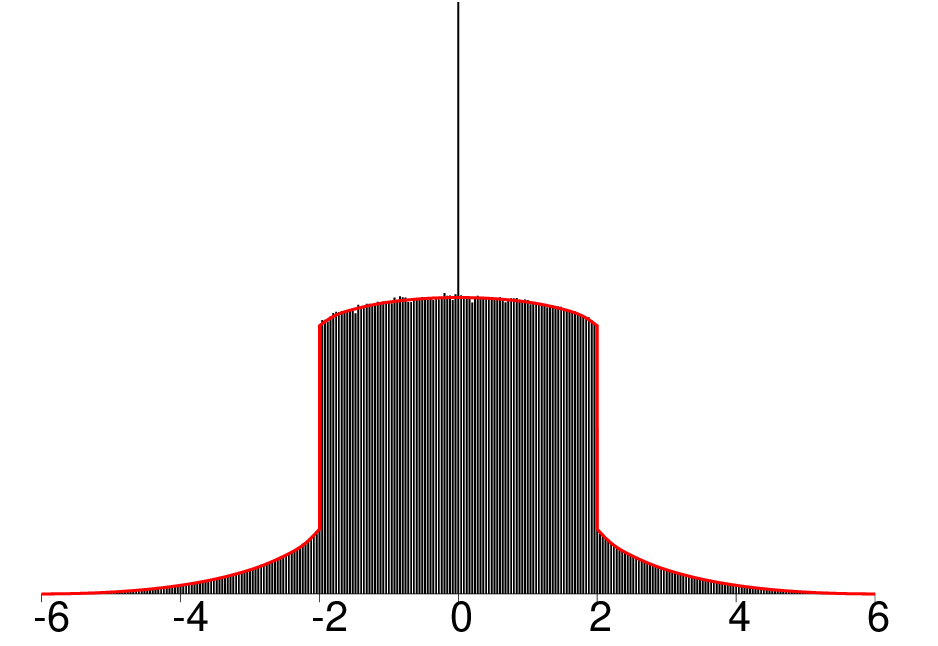}
\end{center}\vspace{-5mm}
\caption{Trace distributions for Examples \ref{ex5} and~\ref{ex6}}
\end{figure}\vspace{-2mm}

In the figure above, each spike is of
mass~$1/2$.
As~before, the~red lines show the densities of the theoretical trace distributions. They~are obtained as the superpositions of the distributions for the two components
$[\SO_3(\bbR)]^2$
and
$[\O_3^-(\bbR)]^2$.
The~first is constructed as explained in Section~\ref{Lie}, the second one by mirroring on the
\mbox{$y$-axis}.

\subsubsection*{An example with CM by an endomorphism field of degree six}

\begin{ex}
\label{ex7}
Let~$X_7'$
be the double cover
of~$\Pb^2_\bbQ$,
given by
$$w^2 = xyz (x^3 - 3x^2z - 3xy^2 - 3xyz + y^3 + 9y^2z + 6yz^2 + z^3)$$
and
$X_7$
the
$K3$~surface
obtained as the minimal desingularisation
of~$X_7'$.

\begin{abc}
\item
Then the geometric Picard rank
of~$X_7$
is
$16$.
\item
The endomorphism field
of~$X_7$
contains~$\bbQ(\sqrt{-1})$.
\end{abc}\smallskip

\noindent
{\bf Proof.}
a)
One~has a lower bound
of~$16$,
as the ramification locus has 15 singular~points. An~upper bound
of~$16$
is provided by the reduction
modulo~$5$,
which is of geometric Picard
rank~$16$.\smallskip

\noindent
b)
The automorphism
of~$\Pb^2_{\overline\bbQ}$,
given by the matrix
$$
\left(
\begin{array}{rrr}
-g\!+\!1 & g^2\!+\!g\!-\!1 & g^2 \\
       0 &               0 &  -g \\
       1 &               1 &   1
\end{array}
\right) \,,
$$
for
$g := \zeta_9+\zeta_9^{-1}-1$,
transforms
$X_7$
into a fibre of the family
$q\colon \calX \to B$,
considered in Theorem~\ref{Qi_deg2}. The assertion follows, as the endomorphism field does not shrink under specialisation.
\eop
\end{ex}

There is strong evidence that
$X_7$
has complex multiplication by the endomorphism field
$E = \bbQ(\zeta_9 + \zeta_9^{-1},\sqrt{-1})$,
which is abelian of degree six. The evidence has been described in \cite[last subsection]{EJ16}. Note~that
$\smash{X_7 = V^{(-1,\mu_9)}}$
in the notation of \cite[Conjectures~5.2]{EJ16}. Thus, conjecturally,
$\smash{(\ST^2(X_7))^0 \cong [\U_1]^3}$.

The~observation that
$\smash{\Tr(\varrho^2_{X_7,l,\tr}(\Frob_p)) = 0}$
for all primes
$p \not\equiv \pm1 \pmod{36}$
has meanwhile been extended to
$p < 10^8$.
If~one knew this unconditionally then Example~\ref{forbidden_elements}.b) would show that
$\ST^2_\tr(X_7)/(\ST^2(X_7))^0 \cong \bbZ/6\bbZ$.
Note~that the maximal totally real subfield
$E_0 = \bbQ(\zeta_9 + \zeta_9^{-1})$
of the conjectural endomorphism field is cyclic of
degree~$3$.\vspace{-4.5mm}

\begin{figure}[H]
\begin{center}
\includegraphics[scale=0.47]{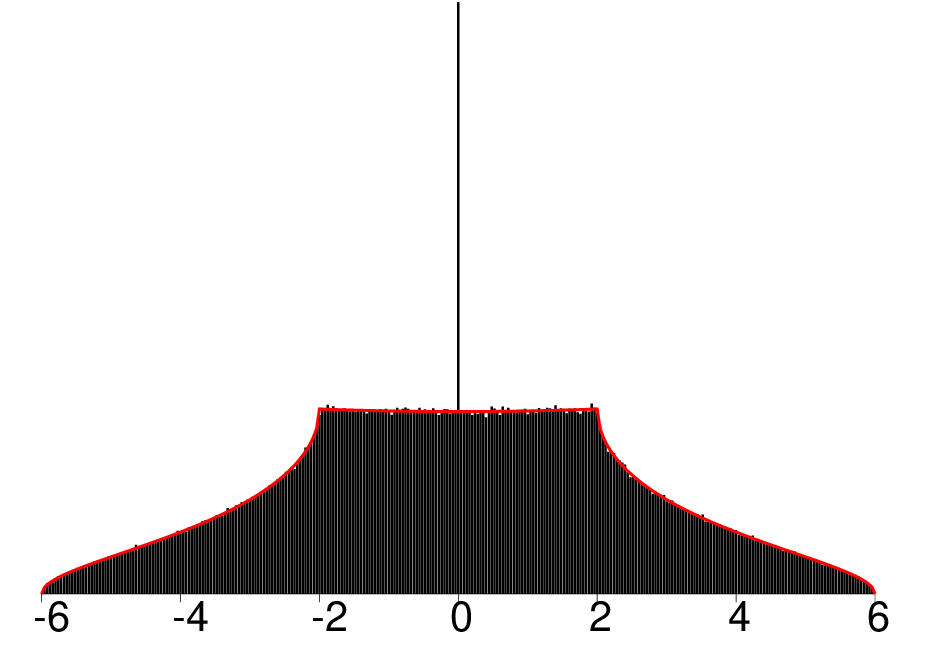}
\end{center}\vspace{-5mm}
\caption{Trace distribution for Example~\ref{ex7}}\label{CM_dim_1_plot}
\end{figure}\vspace{-3mm}

In the figure above, the spike is of
mass~$5/6$.

\subsubsection*{Conclusion}

For each of the seven
$K3$~surfaces
in the sample, we see a strong coincidence in the data that supports the Sato--Tate conjecture.

\subsubsection*{The order of convergence}

In~Example \ref{ex7}, up to
$10^8$,
exactly
$960\,272$
of the
$5\,761\,455$
primes do not contribute to the~spike. Only~these are to be considered. Then,~among the 300 subintervals, the largest discrepancy between the experimental count of Frobenius traces and the theoretical prediction occurs in the subinterval ranging from
$1.72$
to~$1.76$.
Here,
$5538.39$
Frobenius traces are to be expected, but only
$5341$
are~found, a relative error of
roughly~$3.5\%$.

\begin{table}[H]
\begin{center}
{\tiny
\begin{tabular}{|r||l|c|c|}\hline
   $N := $ \hspace{0.25cm} & \quad Largest     & Maximal \# & Largest \\
  \#primes \hspace{0.25cm} & \quad discrepancy & of traces  & discrepancy \\
           \hspace{0.25cm} & \quad             & expected   & $/\sqrt{\text{\#expected}}$ \\\hline\hline
      $32$ \hspace{0.25cm} & \quad $\phantom{-00}1.918$ &  $\phantom{000}0.19$ &  $\phantom{-}4.46$ \\\hline
      $64$ \hspace{0.25cm} & \quad $\phantom{-00}1.836$ &  $\phantom{000}0.37$ &  $\phantom{-}3.02$ \\\hline
     $128$ \hspace{0.25cm} & \quad $\phantom{-00}3.261$ &  $\phantom{000}0.74$ &  $\phantom{-}3.79$ \\\hline
     $256$ \hspace{0.25cm} & \quad $\phantom{-00}3.804$ &  $\phantom{000}1.48$ &  $\phantom{-}3.13$ \\\hline
     $512$ \hspace{0.25cm} & \quad $\phantom{-00}5.092$ &  $\phantom{000}2.96$ &  $\phantom{-}2.96$ \\\hline
  $1\,024$ \hspace{0.25cm} & \quad $\phantom{-00}6.156$ &  $\phantom{000}5.93$ &  $\phantom{-}2.53$ \\\hline
  $2\,048$ \hspace{0.25cm} & \quad $\phantom{-00}9.285$ &  $\phantom{00}11.85$ &  $\phantom{-}2.70$ \\\hline
  $4\,096$ \hspace{0.25cm} & \quad $\phantom{-0}14.364$ &  $\phantom{00}23.70$ &  $\phantom{-}2.95$ \\\hline
  $8\,192$ \hspace{0.25cm} & \quad $\phantom{-0}23.728$ &  $\phantom{00}47.40$ &  $\phantom{-}3.44$ \\\hline
 $16\,384$ \hspace{0.25cm} & \quad $\phantom{-0}30.457$ &  $\phantom{00}94.80$ &  $\phantom{-}3.13$ \\\hline
 $32\,768$ \hspace{0.25cm} & \quad $\phantom{-0}41.965$ &  $\phantom{0}189.60$ &  $\phantom{-}3.05$ \\\hline
 $65\,536$ \hspace{0.25cm} & \quad $\phantom{-0}53.930$ &  $\phantom{0}379.21$ &  $\phantom{-}2.77$ \\\hline
$131\,072$ \hspace{0.25cm} & \quad $\phantom{-0}75.102$ &  $\phantom{0}758.42$ &  $\phantom{-}2.73$ \\\hline
$262\,144$ \hspace{0.25cm} & \quad $-129.287$ & $1516.84$ & $-3.32$ \\\hline
$524\,288$ \hspace{0.25cm} & \quad $-146.727$ & $3033.68$ & $-2.66$ \\\hline
$960\,272$ \hspace{0.25cm} & \quad $-197.392$ & $5556.40$ & $-2.65$ \\\hline
\end{tabular}
}
\end{center}\vspace{-2.5mm}

\caption{Discrepancies between numbers of Frobenius traces}
\end{table}\vspace{-3.5mm}

The~maximal number of traces to be expected in a subinterval is
$5556.40$.
This number does not occur near
$0.00$,
but for the subintervals
$[-2.00, -1.96]$
and~$[1.96, 2.00]$.
One~calculates that
$(5341-5538.39)/\sqrt{5556.40} \approx -2.65$.

Doing~the same for only the first
$2^k$
good primes, not contributing to the spike, for
$k=3,\ldots,19$,
the data were obtained that are presented in Table~4 above. It~seems that the values in the column to the right remain within a bounded range around~zero.

As~the maximal number of traces expected among the subintervals is proportional to the number of
primes~$N$,
this suggests that the largest discrepancy is proportional
to~$\sqrt{N}$.
Consequently,~the
\mbox{$L^0$-distance}
between the experimental and theoretical density functions is proportional
to~$\smash{\frac1{\sqrt{N}}}$,
which means that convergence is of
order~$\smash{\frac12}$.

\subsubsection*{The order of convergence--A second experiment}

Experts in Statistics advise to consider the
\mbox{$L^0$-distance}
between the experimental and theoretical distribution functions, instead of the densities. Again, the values in the column to the right seem to fluctuate within a bounded range around~zero. Which~would indeed show convergence of
order~$\smash{\frac12}$.\vspace{-2mm}

\begin{table}[H]
\begin{center}
{\tiny
\begin{tabular}{|r||l|c|}\hline
           \hspace{0.25cm} & \quad $L^0$-distance & \\
   $N := $ \hspace{0.25cm} & \quad between        & multiplied \\
  \#primes \hspace{0.25cm} & \quad distribution   & by $\sqrt{N}$ \\
           \hspace{0.25cm} & \quad functions      & \\\hline\hline
       $8$ \hspace{0.25cm} & \quad $-0.175$ & $-0.496$ \\\hline
      $16$ \hspace{0.25cm} & \quad $-0.262$ & $-1.049$ \\\hline
      $32$ \hspace{0.25cm} & \quad $-0.155$ & $-0.877$ \\\hline
      $64$ \hspace{0.25cm} & \quad $-0.077\,5$ & $-0.620$ \\\hline
     $128$ \hspace{0.25cm} & \quad $-0.047\,2$ & $-0.534$ \\\hline
     $256$ \hspace{0.25cm} & \quad $\phantom{-}0.057\,2$ & $\phantom{-}0.916$ \\\hline
     $512$ \hspace{0.25cm} & \quad $\phantom{-}0.034\,1$ & $\phantom{-}0.772$ \\\hline
  $1\,024$ \hspace{0.25cm} & \quad $-0.018\,8$ & $-0.601$ \\\hline
  $2\,048$ \hspace{0.25cm} & \quad $-0.013\,1$ & $-0.595$ \\\hline
  $4\,096$ \hspace{0.25cm} & \quad $-0.007\,40$ & $-0.473$ \\\hline
  $8\,192$ \hspace{0.25cm} & \quad $-0.007\,05$ & $-0.638$ \\\hline
 $16\,384$ \hspace{0.25cm} & \quad $\phantom{-}0.004\,12$ & $\phantom{-}0.527$ \\\hline
 $32\,768$ \hspace{0.25cm} & \quad $-0.003\,55$ & $-0.642$ \\\hline
 $65\,536$ \hspace{0.25cm} & \quad $\phantom{-}0.003\,06$ & $\phantom{-}0.783$ \\\hline
$131\,072$ \hspace{0.25cm} & \quad $\phantom{-}0.001\,70$ & $\phantom{-}0.615$ \\\hline
$262\,144$ \hspace{0.25cm} & \quad $-0.001\,50$ & $-0.770$ \\\hline
$524\,288$ \hspace{0.25cm} & \quad $-0.001\,21$ & $-0.874$ \\\hline
$960\,272$ \hspace{0.25cm} & \quad $-0.000\,521$ & $-0.511$ \\\hline
\end{tabular}
}
\end{center}\vspace{-2.5mm}

\caption{$L^0$-distance
between experimental and theoretical distribution~functions}
\end{table}\vspace{-2.2mm}

For~the other six example surfaces, we made the analogous experiments. We~do not think that it is useful to present the corresponding raw data in this article. In~fact, all the surfaces examined show qualitatively the same~behaviour.

\subsubsection*{The Lang--Trotter conjecture}\leavevmode\medskip

\noindent
Concerning the traces of the Frobenii
$\Frob_p$
on
$\smash{H^2_\et(X_{7,\overline\bbQ}, \bbQ_l)}$,
the statistics for
$p<10^8$
is as~follows. There~are
$922\,644$
distinct integers occurring as a~trace. One~of them
is~$0$,
which comes up roughly
$\smash{\frac56}$
of the~time. Except for this, there are
$886\,932$
integers that occur only once,
$33\,922$
integers that occur exactly twice,
$1\,667$
integers that occur exactly three times,
$115$
integers that occur exactly four times, and
$7$
integers that occur exactly five times as a~trace. No~integer occurs more than five~times. Comparing~this with
$\log \log 10^8 \approx 2.913$,
there is certainly no contradiction with the Lang--Trotter conjecture to be seen from our~data.
Once~again, the other examples show qualitatively the same~behaviour.

\appendix
\section{A family that is acted upon by
$\bbQ(\sqrt{-1})$}

\begin{theo}
\label{Qi_deg2}
Let\/~$B \subset \Pb^2_\bbQ \!\times\! \Pb^2_\bbQ$
be the closed subscheme given by the equations\/
$a_1 b_3 + a_2 b_1 - 2a_3 b_1 = 0$
and\/
$a_1 b_2 + a_2 b_3 - 2a_3 b_2 = 0$,
and let, moreover,
$q\colon \calX' \to B$
be the family of double covers
of\/~$\Pb^2$
given~by
$$w^2 = l_1 \cdots l_6 \,,$$
for\/
$l_1, \ldots, l_6$
the linear forms\/
$l_1 := x$,
$l_2 := y$,
$l_3 := z$,
$l_4 := x + y + z$,
$l_5 := a_1 x + a_2 y + a_3 z$,
and\/~$l_6 := b_1 x + b_2 y + b_3 z$.

\begin{abc}
\item
Then the generic fibre
$\calX'_\eta$
is normal surface, the minimal desingularisation of which is a\/
$K3$~surface\/
$\calX_\eta$
of geometric Picard
rank\/~$16$.
\item
The endomorphism field
of\/~$\calX_\eta$
is\/~$\bbQ(\sqrt{-1})$.
\end{abc}\smallskip

\noindent
{\bf Proof.}
{\em
a)
The singularities
of~$\calX'_\eta$
are caused by those of the ramification curve, and are therefore isolated and of type~$A_1$.
The~surface
$\calX_\eta$
is a
$K3$~surface
as the ramification curve is of degree six. For~the geometric Picard rank, one has a lower bound
of~$16$,
as the ramification locus has 15 singular~points. An~upper bound
of~$16$
is provided by the specialisation
to~$X_4$,
cf.\ Example~\ref{ex4}.c).\smallskip

\noindent
b)
The specialisation
to~$X_4$
is known to have an endomorphism field of
degree~$\leq\! 2$,
cf.~Example~\ref{ex4}.d). As~the endomorphism field does not shrink under specialisation \cite[Corollary~4.6]{EJ20a}, it is sufficient to show that the endomorphism field
of~$\smash{\calX_\eta}$
contains~$\smash{\bbQ(\sqrt{-1})}$.

For this, we blow up
$\smash{\Pb^2_{k(B)}}$
in the seven points
$\Vb(l_i,l_j)$,
for
$\{i,j\} = \{1,2\}$,
$\{1,4\}$,
$\{1,6\}$,
$\{2,4\}$,
$\{2,6\}$,
$\{4,6\}$,
and~$\{3,5\}$.
Since~no four of these points are collinear, the result is a weak del Pezzo surface~$S$
of
degree~$2$~\cite[Corollary~8.1.24]{Do}.
The~linear system of the cubic forms vanishing in the seven blown-up points defines a birational morphism
$S \to S'$
to a singular model
$S'\colon W^2 = Q(X,Y,Z)$.
Here,~$Q$
defines a plane quartic having only simple singularities~\cite[Theorem~8.3.2.(iv)]{Do}. A~calculation shows that, in our particular situation, the quartic
$\Vb(Q)$
splits over
$k(B)$
into the union of two conics,
$Q = Q_1Q_2$.

The double cover
$\smash{\calX'_\eta}$
of~$\smash{\Pb^2_{k(B)}}$
goes over, under blowing up, into a double cover
of~$S$,
and therefore also into one
of~$S'$.
The~special choice
of~$B$
makes sure that the ramification locus
$\Vb(l_1 \cdots l_6)$
is mapped
to~$\Vb(Q)$.
In~other words, a linear algebra calculation over the function field
$k(B)$
shows that three linearly independent cubic forms vanishing in the seven blown-up points, together with the
coordinate~$w$,
fulfil exactly one quartic relation, which is of the kind
$W^4 = Q(X,Y,Z)$.
I.e.,
$\calX_\eta$
has a singular
model~$\calX''_\eta$
of degree four, which is given by the~equation
$$W^4 = Q_1(X,Y,Z)Q_2(X,Y,Z)\,.$$
There~is an automorphism
of~$\calX''_\eta$,
given by
$I\colon (W\!:\!X\!:\!Y\!:\!Z) \mapsto (iW\!:\!X\!:\!Y\!:\!Z)$,
cf.\ \cite[Example~6.17]{EJ20b}. We~claim that the operation
of~$I$
on~$H_\tr$
gives rise to complex~multiplication. For~this, it has to be shown
that~$J := I \!\circ\! I$
acts
on~$H_\tr$
as the multiplication
by~$(-1)$.
Let~us note that
$H_\tr$
is a simple
\mbox{$\bbQ$-Hodge}
structure~\cite[Theorem~1.6.a)]{Za}, so it suffices to exclude multiplication
by~$(+1)$.

For~this, observe that
$\smash{\calX''_\eta}$
has four singularities, each of which is of
type~$A_3$.
Hence,
$\smash{\dim H^2(\calX''_\eta(\bbC), \bbQ) = 10}$.
The~fixed point set
of~$J$
is the union of two conics, which has topological Euler
characteristic~$0$.
Therefore, the Lefschetz trace formula \cite[Theorem 8.5]{Ed} shows that
$\smash{\Tr(J|_{H^2(\calX''_\eta(\bbC), \bbQ)}) = -2}$.
In other words,
$\smash{J|_{H^2(\calX''_\eta(\bbC), \bbQ)}}$
has the eigenvalue
$(+1)$
with
multiplicity~$4$,
while the eigenvalue
$(-1)$
occurs with
multiplicity~$6$.
In~particular,
$\smash{H_\tr}$,
which is of dimension six, cannot be contained in the
\mbox{$(+1)$-eigenspace},
which completes the~proof.
}
\eop
\end{theo}

\section{A table}

\begin{table}[H]\vspace{-2mm}
\begin{tabular}{|c||c|c|c|c|c|c|c|}\hline
$i$                 & 1                  & 2              & 3              & 4               & 5              & 6                & 7               \\\hline\hline
Bad primes of $X_i$ & $2,3,5,7,$         & $2,3,5,7$      & $2,3,5$        & $2,3,5,$        & $2,5$          & $2,3,5,7,$       & $2,3$           \\
                    & $11,13,29$         &                &                & $7,11$          &                & $11,13,17,47$    &                 \\\hline
Jump character      &                    &                &                &                 &                &                  &                 \\
$\Delta_\tr \in \bbQ^*/\bbQ^{*2}$ of $X_i$ 
                    & $-\overline{6006}$ & $\overline{1}$ & $\overline{1}$ & $-\overline{1}$ & $\overline{1}$ & $(\overline{3})$ & $-\overline{1}$ \\\hline
\end{tabular}%\smallskip

\caption{Bad primes and jump character for the surfaces in the sample}\label{badpr_jump}
\end{table}

\setlength\parindent{0mm}
\end{document}